\documentclass[a4paper,11pt,reqno,twoside]{amsart}
\usepackage{riro_hal}
\begin{document} 
%..........................................
\begin{merci}
This work began while both authors shared the hospitality of Centre de %
Recherches Math\'ematiques (Montr\'eal) during the %
theme year "Analysis in Number Theory" (first semester of 2006). We would like %
to thank C.~David, H.~Darmon and %
A.~Granville for their invitation. This work was essentially completed in %
september 2006 in CIRM (Luminy) at %
the occasion of J.-M.~Deshouillers' sixtieth birthday. %
We would like to wish him the best. 
\end{merci}
%..........................................
\section{Introduction and statement of the results}

\subsection{Description of the families of $L$-functions studied}
The purpose of this paper is to compute various statistics associated to low-lying zeros of several families %
of symmetric power $L$-functions in the level aspect. %
First of all, we give a short description of these families. %
To any primitive holomorphic cusp form $f$ of prime level $q$ %
and even weight%
\footnote{In this paper, the weight $\kappa$ is a \emph{fixed} even integer and the level $q$ goes %
to infinity among the prime numbers.} $\kappa\geq 2$ %
(see \S~\ref{sec_autoback} for the automorphic background) %
say $f\in\prim{\kappa}{q}$, one can associate its $r$-th symmetric power $L$-function denoted by $L(\sym^rf,s)$ for any %
integer $r\geq 1$. It is given by an explicit absolutely convergent Euler product of degree $r+1$ on $\Re{s}>1$ %
(see \S~\ref{sec_sympow}). The completed $L$-function is defined by
\begin{equation*}%
\Lambda(\sym^rf,s)\coloneqq\left(q^{r}\right)^{s/2}L_\infty(\sym^rf,s)L(\sym^rf,s)%
\end{equation*}
where $L_\infty(\sym^rf,s)$ is a product of $r+1$ explicit $\fGamma_{\R}$-factors (see \S~\ref{sec_sympow}) %
and $q^r$ is the arithmetic conductor. We will need some control on the analytic behaviour of this %
function. Unfortunately, such information is not currently %
known in all generality. We sum up our main assumption in the following statement. \label{hypohypo}
\begin{hypothesis}
The function $\Lambda\left(\sym^rf,s\right)$ is a \emph{completed $L$-function} in the
sense that it satisfies the 
following \emph{nice} analytic properties:
\begin{itemize}
\item
it can be extended to an holomorphic function of order $1$ on $\C$,
\item
it satisfies a functional equation of the shape
\[%
\Lambda(\sym^rf,s)=\epsilon\left(\sym^rf\right)\Lambda(\sym^rf,1-s)
\]
where the sign $\epsilon\left(\sym^rf\right)=\pm 1$ of the functional equation is given by
\begin{equation}\label{valueofsign}%
\epsilon\left(\sym^rf\right)\coloneqq%
\begin{cases}
+1 & \text{if $r$ is even},\\
\epsilon_f(q)\times\epsilon(\kappa,r) & \text{otherwise}
\end{cases}
\end{equation}
with
\[%
\epsilon(\kappa,r)\coloneqq i^{\left(\frac{r+1}{2}\right)^2(\kappa-1)+\frac{r+1}{2}}=%
\begin{cases}
i^{\kappa} & \text{if $\;r\equiv 1\pmod{8}$,} \\
-1 & \text{if $\;r\equiv 3\pmod{8}$}, \\
-i^{\kappa} & \text{if $\;r\equiv 5\pmod{8}$}, \\
+1 & \text{if $r\;\equiv 7\pmod{8}$}
\end{cases}
\]
and $\epsilon_f(q)=\pm 1$ is defined in \eqref{eq_signe} and only depends on $f$ and $q$.
\end{itemize}
\end{hypothesis}
\begin{remint}
Hypothesis $\Nice(r,f)$ is known for $r=1$ %
(E.~Hecke \cite{MR1513069,MR1513122,MR1513142}), $r=2$ thanks to the work of S.~Gelbart and H.~Jacquet \cite{GeJa} %
and $r=3,4$ from the works of H.~Kim and F.~Shahidi \cite{KiSh1,KiSh2,Ki}.
\end{remint}
We aim at studying the low-lying zeros for the family of $L$-functions given by
\begin{equation*}
\mathcal{F}_r\coloneqq%
\bigcup_{\text{$q$ prime}}\left\{L(\sym^rf,s), f\in\prim{\kappa}{q}\right\}
\end{equation*}
for any integer $r\geq 1$. Note that when $r$ is even, the sign of the functional equation of any $L(\sym^rf,s)$ %
is constant of value $+1$ but when $r$ is odd, this is definitely not the case. As a consequence, it is very natural %
to understand the low-lying zeros for the subfamilies given by
\begin{equation*}
\mathcal{F}_r^{\epsilon}\coloneqq%
\bigcup_{\text{$q$ prime}}\left\{L(\sym^rf,s), \;f\in\prim{\kappa}{q}, \epsilon\left(\sym^rf\right)=\epsilon\right\}
\end{equation*}
for any odd integer $r\geq 1$ and for $\epsilon=\pm 1$.

\subsection{Symmetry type of these families}

One of the purpose of this work is to determine the symmetry type of the families $\mathcal{F}_r$ and $\mathcal{F}_r^{\epsilon}$ %
for $\epsilon=\pm 1$ and for any integer $r\geq 1$ (see \S~\ref{sec_densres} for the background on symmetry types). %
The following theorem is a quick summary of the symmetry types obtained. 
\begin{theoint}\label{thm_A}%
Let $r\geq 1$ be any integer and $\epsilon=\pm 1$. We assume that hypothesis $\Nice(r,f)$ holds for any prime number $q$ and any %
primitive holomorphic cusp form of level $q$ and even weight $\kappa\geq 2$. The symmetry group $G(\mathcal{F}_r)$ of $\mathcal{F}_r$ is %
given by
\begin{equation*}
G(\mathcal{F}_r)=\begin{cases}%
Sp & \text{ if $r$ is even,}\\%
O & \text{ otherwise.}%
\end{cases}%
\end{equation*}
If $r$ is odd then the symmetry group $G(\mathcal{F}_r^\epsilon)$ of $\mathcal{F}_r^\epsilon$ is given by %
\begin{equation*}
G(\mathcal{F}_r^\epsilon)=\begin{cases}%
SO(\mathrm{even}) & \text{ if $\epsilon=+1$,}\\%
SO(\mathrm{odd}) & \text{ otherwise.}%
\end{cases}
\end{equation*}
\end{theoint}
\begin{remint}
\label{remark2}
It follows in particular from the value of $\epsilon\left(\sym^rf\right)$ given in \eqref{valueofsign} that, if $r$ is even, %
then $\sym^rf$ has not the same symmetry type %
than $f$ and, if $r$ is odd, then $f$ and $\sym^rf$ have the same symmetry type if and only if
\begin{align*}
r \equiv 1 \pmod{8}\; \text{ and }\; \kappa\equiv 0\pmod{4}%
\shortintertext{or}%
r \equiv 5 \pmod{8}\; \text{ and }\; \kappa\equiv 2\pmod{4}%
\shortintertext{or}%
r \equiv 7 \pmod{8}.%
\end{align*}
\end{remint}
\begin{remint}
Note that we do not assume any Generalised Riemann Hypothesis for the symmetric power $L$-functions.
\end{remint}
In order to prove theorem~\ref{thm_A}, we compute either the (signed) asymptotic expectation of the one-level density or %
the (signed) asymptotic expectation of the two-level density. The results are given in the next two sections in %
which $\epsilon=\pm 1$, $\nu$ will always be a positive real number, $\Phi, \Phi_1$ and $\Phi_2$ will always stand for %
even Schwartz functions whose Fourier transforms $\widehat{\Phi}, \widehat{\Phi_1}$ and $\widehat{\Phi_2}$ are compactly supported in $[-\nu,+\nu]$ and $f$ will always be a %
primitive holomorphic cusp form of prime level $q$ and even weight $\kappa\geq 2$ for which hypothesis $\Nice(r,f)$ holds. %
We refer to \S~\ref{sec_proba} for the probabilistic background.

\subsubsection{(Signed) asymptotic expectation of the one-level density}

The \emph{one-level density} (relatively to $\Phi$) of $\sym^rf$ is defined by
\begin{equation*}
D_{1,q}[\Phi;r](f)%
\coloneqq%
\sum_{\rho,\;\Lambda(\sym^rf,\rho)=0}%
\Phi\left(\frac{\log{\left(q^r\right)}}{2i\pi}\left(%
\Re{\rho}-\frac{1}{2}+i\Im{\rho}%
\right)\right)%
\end{equation*}
where the sum is over the non-trivial zeros $\rho$ of $L(\sym^rf,s)$ with multiplicities. The \emph{asymptotic expectation} %
of the one-level density is by definition %
\begin{equation*}
\lim_{\substack{q\;\text{prime} \\ q\to+\infty}}%
\smashoperator[r]{%
\sum_{f\in\prim{\kappa}{q}}%
}\omega_q(f)D_{1,q}[\Phi;r](f)
\end{equation*}
where $\omega_q(f)$ is the harmonic weight defined in \eqref{eq_facomeg} and similarly the \emph{signed asymptotic expectation} of %
the one-level density is by definition %
\begin{equation*}
\lim_{\substack{q\;\text{prime} \\ q\to+\infty}}%
2%
\smashoperator[r]{%
\sum_{\substack{f\in\prim{\kappa}{q} \\ \epsilon\left(\sym^rf\right)=\epsilon}}%
}
\omega_q(f)D_{1,q}[\Phi;r](f)
\end{equation*}
when $r$ is odd.
\begin{theoint}\label{thm_B}%
Let $r\geq 1$ be any integer and $\epsilon=\pm 1$. We assume that hypothesis $\Nice(r,f)$ holds for %
any prime number $q$ and any primitive holomorphic cusp form of level $q$ and even weight $\kappa\geq 2$ and also %
that $\theta$ is admissible (see hypothesis $\Hy_2(\theta)$ %
page~\pageref{hyp_RPS}). Let%
\[%
\nu_{1,\mathrm{max}}(r,\kappa,\theta)\coloneqq%
\left(1-\frac{1}{2(\kappa-2\theta)}\right)\frac{2}{r^2}.%
\]
If $\nu<\nu_{1,\mathrm{max}}(r,\kappa,\theta)$ then the asymptotic expectation of the one-level density is
\begin{equation*}
\widehat{\Phi}(0)+\frac{(-1)^{r+1}}{2}\Phi(0).
\end{equation*}
Let
\[%
\nu_{1,\mathrm{max}}^\epsilon(r,\kappa,\theta)\coloneqq%
\inf{\left(\nu_{1,\mathrm{max}}(r,\kappa,\theta),\frac{3}{r(r+2)}\right).}%
\]
If $r$ is odd and $\nu<\nu_{1,\mathrm{max}}^\epsilon(r,\kappa,\theta)$ then the signed asymptotic expectation %
of the one-level density is %
\begin{equation*}%
\widehat{\Phi}(0)+\frac{(-1)^{r+1}}{2}\Phi(0).%
\end{equation*}
\end{theoint}
\begin{remint}
\label{remark4}
The first part of Theorem~\ref{thm_B} reveals that the symmetry type of $\mathcal{F}_r$ is
\[%
G(\mathcal{F}_r)=%
\begin{cases}
Sp & \text{if $r$ is even,}\\%
O & \text{if $r=1$,}\\%
SO(\mathrm{even}) \text{ or } O \text{ or } SO(\mathrm{odd}) & \text{if $r\geq 3$ is odd.}%
\end{cases}%
\]
We cannot decide between the three orthogonal groups when $r\geq 3$ is odd since in this %
case $\nu_{1,\mathrm{max}}(r,\kappa,\theta)<1$ but the computation of the two-level densities will %
enable us to decide. Note also that we go beyond the support $[-1,1]$ when $r=1$ %
as Iwaniec, Luo \& Sarnak \cite{IwLuSa} (Theorem 1.1) but without doing any subtle % 
arithmetic analysis of Kloosterman sums. Also, A.~G{\"u}loglu % 
in \cite[Theorem 1.2]{Gu} established some density result for %
the same family of $L$-functions but when the weight $\kappa$ goes to infinity and the level $q$ is fixed. It %
turns out that we recover the same %
constraint on $\nu$ when $r$ is even but we get a better result when $r$ is %
odd. This can be explained by the fact %
that the analytic conductor of any $L(\sym^rf,s)$ with $f$ in $\prim{\kappa}{q}$ which is of size %
\[%
q^r\times%
\begin{cases}%
\kappa^{r} & \text{if $r$ is even}\\%
\kappa^{r+1} & \text{otherwise}%
\end{cases}%
\]
is slightly larger in his case than in ours when $r$ is odd.
\end{remint}
\begin{remint}
\label{remark5}
The second part of Theorem~\ref{thm_B} reveals that if $r$ is odd and $\epsilon=\pm 1$ then the symmetry type of $\mathcal{F}_r^\epsilon$ is %
\[%
G(\mathcal{F}_r^{\epsilon})=%
SO(\mathrm{even}) \text{ or } O \text{ or } SO(\mathrm{odd}).%
\]
Here $\nu$ is always strictly smaller than one and we are not able to recover the result of \cite[Theorem 1.1]{IwLuSa} %
without doing some arithmetic on Kloosterman sums.
\end{remint}

\subsubsection{Sketch of the proof}

We give here a sketch of the proof of the first part of Theorem \ref{thm_B} namely we briefly explain how to determine the asymptotic expectation of the one-level density assuming that hypothesis $\Nice(r,f)$ holds for any prime number $q$ and any primitive holomorphic cusp form of level $q$ and even weight $\kappa\geq 2$ and also that $\theta$ is admissible. The first step consists in transforming the sum over the zeros of $\Lambda(\sym^rf,s)$ which occurs in $D_{1,q}[\Phi;r](f)$ into a sum over primes. This is done \emph{via} some Riemann's explicit formula for symmetric power $L$-functions stated in Proposition \ref{explicit} which leads to
\begin{equation*}
D_{1,q}[\Phi;r](f)=\widehat{\Phi}(0)+\frac{(-1)^{r+1}}{2}\Phi(0)+P_q^1[\Phi;r](f)+\sum_{m=0}^{r-1}(-1)^mP_q^2[\Phi;r,m](f)+o(1)
\end{equation*}
where
\begin{equation}\label{eq_cneser}
P_{q}^1[\Phi;r](f)\coloneqq-\frac{2}{\log{\left(q^r\right)}}\sum_{\substack{p\in\prem \\ p\nmid q}}\lambda_{f}\left(p^r\right)\frac{\log{p}}{\sqrt{p}}\widehat{\Phi}\left(\frac{\log{p}}{\log{\left(q^r\right)}}\right).
\end{equation}
The terms $P_q^2[\Phi;r,m](f)$ are also sums over primes which look like $P_{q}^1[\Phi;r](f)$ but can be forgotten in first approximation since they can be thought as sums over squares of primes which are easier to deal with. The second step consists in averaging over all the $f$ in $\prim{\kappa}{q}$. While doing this, the asymptotic expectation of the one-level density
\begin{equation*}
\widehat{\Phi}(0)+\frac{(-1)^{r+1}}{2}\Phi(0)
\end{equation*}
naturally appears and we need to show that
\begin{equation*}
-\frac{2}{\log{\left(q^r\right)}}\sum_{\substack{p\in\prem \\ p\nmid q}}\left(\sum_{f\in\prim{\kappa}{q}}\omega_q(f)\lambda_{f}\left(p^r\right)\right)\frac{\log{p}}{\sqrt{p}}\widehat{\Phi}\left(\frac{\log{p}}{\log{\left(q^r\right)}}\right)
\end{equation*}
is a remainder term provided that the support $\nu$ of $\Phi$ is small enough. We apply some suitable trace formula given in Proposition \ref{iwlusatr} in order to express the previous average of Hecke eigenvalues. We cannot directly apply Peterson's trace formula since there may be some old forms of level $q$ especially when the weight $\kappa$ is large. Nevertheless, these old forms are automatically of level $1$ since $q$ is prime and their contribution remains negligible. So, we have to bound
\begin{equation*}
-\frac{4\pi i^\kappa}{\log{\left(q^r\right)}}\sum_{\substack{p\in\prem \\ p\nmid q}}\sum_{\substack{c\geq 1 \\
q\mid c}}\frac{S(1,p^r;c)}{c}J_{\kappa-1}\left(\frac{4\pi\sqrt{p^r}}{c}\right)\frac{\log{p}}{\sqrt{p}}\widehat{\Phi}\left(\frac{\log{p}}{\log{\left(q^r\right)}}\right)
\end{equation*}
where $S(1,p^r;c)$ is a Kloosterman sum and which can be written as
\begin{equation*}
-\frac{4\pi i^\kappa}{\log{\left(q^r\right)}}\sum_{\substack{c\geq 1 \\
q\mid c}}\sum_{m\geq 1}a_m\frac{S(1,m;c)}{c}g(m;c)
\end{equation*}
where
\begin{equation*}
a_m\coloneqq\un_{[1,q^{r^2\nu}]}(m)\frac{\log{m}}{rm^{1/(2r)}}\times\begin{cases}
1 & \text{if $m=p^r$ for some prime $p\neq q$}, \\
0 & \text{otherwise}
\end{cases}
\end{equation*}
and
\begin{equation*}
g(m;c)\coloneqq J_{\kappa-1}\left(\frac{4\pi\sqrt{m}}{c}\right)\widehat{\Phi}\left(\frac{\log{m}}{r\log{\left(q^r\right)}}\right).
\end{equation*}
We apply the large sieve inequality for Kloosterman sums given in proposition \ref{sieve}. %
It entails that if $\nu\leq 2/r^2$ then such quantity is bounded by
\begin{equation*}
\ll_\epsilon q^{\left(\frac{\kappa-1}{2}-\theta\right)(r^2\nu-2)+\epsilon}+q^{\left(\frac{\kappa}{2}-\theta\right)r^2\nu-\left(\kappa-\frac{1}{2}-2\theta\right)+\epsilon}.
\end{equation*}
This is an admissible error term if $\nu<\nu_{1,\mathrm{max}}(r,\kappa,\theta)$. We focus on the fact that we did any arithmetic analysis of Kloosterman sums to get this result. Of course, the power of spectral theory of automorphic forms is hidden in the large sieve inequalities for Kloosterman sums.

\subsubsection{(Signed) asymptotic expectation of the two-level density}

The \emph{two-level density} of $\sym^rf$ (relatively to $\Phi_1$ and $\Phi_2$) is defined by %
\[%
D_{2,q}[\Phi_1,\Phi_2;r](f)%
\coloneqq%
\sum_{\substack{(j_1,j_2)\in\mathcal{E}(f,r)^2\\ j_1\neq\pm j_2}}%
\Phi_1\left(\widehat{\rho}_{f,r}^{(j_1)}\right)%
\Phi_2\left(\widehat{\rho}_{f,r}^{(j_2)}\right).%
\]
For more precision on the numbering of the zeros, we refer to \S~\ref{sec_explicit}. The \emph{asymptotic expectation} of the two-level density is by definition %
\begin{equation*}
\lim_{\substack{q\;\text{prime} \\ q\to+\infty}}%
\smashoperator[r]{%
\sum_{f\in\prim{\kappa}{q}}%
}%
\omega_q(f)D_{2,q}[\Phi_1,\Phi_2;r](f)%
\end{equation*}
and similarly the \emph{signed asymptotic expectation} of the two-level density is by definition %
\begin{equation*}
\lim_{\substack{q\;\text{prime} \\ q\to+\infty}}%
2%
\smashoperator[r]{%
\sum_{\substack{f\in\prim{\kappa}{q} \\ \epsilon\left(\sym^rf\right)=\epsilon}}%
}
\omega_q(f)D_{2,q}[\Phi_1,\Phi_2;r](f)%
\end{equation*}
when $r$ is odd and $\epsilon=\pm 1$.
\begin{theoint}\label{thm_C}
Let $r\geq 1$ be any integer and $\epsilon=\pm 1$. We assume that hypothesis $\Nice(r,f)$ holds for any prime number $q$ %
and any primitive holomorphic cusp form of level $q$ and even weight $\kappa\geq 2$. %
If $\nu<1/r^2$ then the asymptotic expectation of the two-level density is %
\begin{multline*}
\left[%
\widehat{\Phi_1}(0)%
+%
\frac{(-1)^{r+1}}{2}\Phi_1(0)%
\right]%
\left[%
\widehat{\Phi_2}(0)%
+%
\frac{(-1)^{r+1}}{2}\Phi_2(0)%
\right]%
\\
+2\int_\R\abs{u}\widehat{\Phi_1}(u)\widehat{\Phi_2}(u)\dd u %
-2\widehat{\Phi_1\Phi_2}(0)%
+\left((-1)^{r}+\frac{%
%\delta(2\nmid r)%
\un_{2\N+1}(r)
}{2}\right)\Phi_1(0)\Phi_2(0).%
\end{multline*}
If $r$ is odd and $\nu<1/(2r(r+2))$ then the signed asymptotic expectation of the two-level density is %
\begin{multline*}
\left[%
\widehat{\Phi_1}(0)%
+%
\frac{1}{2}\Phi_1(0)%
\right]%
\left[%
\widehat{\Phi_2}(0)%
+%
\frac{1}{2}\Phi_2(0)%
\right]%
\\%
+2\int_\R\abs{u}\widehat{\Phi_1}(u)\widehat{\Phi_2}(u)\dd u%
-2\widehat{\Phi_1\Phi_2}(0)%
-\Phi_1(0)\Phi_2(0)%
\\%
+%
%\delta(\epsilon=-1)%
\un_{\{-1\}}(\epsilon)
\Phi_1(0)\Phi_2(0).%
\end{multline*}
\end{theoint} 
\begin{remint}
\label{remark6}
We have just seen that the computation of the one-level density already reveals that the symmetry %
type of $\mathcal{F}_r$ is $Sp$ when $r$ is even. The asymptotic expectation of the two-level density %
also coincides with the one of $Sp$ (see \cite[Theorem A.D.2.2]{KaSa} or \cite[Theorem 3.3]{Mil}). %
When $r\geq 3$ is odd, the first part of Theorem~\ref{thm_C} together with a result of Katz \& Sarnak %
(see \cite[Theorem A.D.2.2]{KaSa} or \cite[Theorem 3.2]{Mil}) imply that the symmetry type of $\mathcal{F}_r$ is $O$. % 
\end{remint}
\begin{remint}
The second part of Theorem~\ref{thm_C} and a result of Katz \& Sarnak %
(see \cite[Theorem A.D.2.2]{KaSa} or \cite[Theorem 3.2]{Mil}) imply that the symmetry type of $\mathcal{F}_r^{\epsilon}$ %
is as in Theorem~\ref{thm_A} for any odd integer $r\geq 1$ and $\epsilon=\pm 1$.
\end{remint}
In order to prove Theorem~\ref{thm_C}, we need to determine the \emph{asymptotic variance} of the one-level density which is defined by %
\begin{equation*}
\lim_{\substack{q\;\text{prime} \\ q\to+\infty}}%
\smashoperator[r]{%
\sum_{f\in\prim{\kappa}{q}}%
}%
\omega_q(f)\left(D_{1,q}[\Phi;r](f)-%
\sum_{g\in\prim{\kappa}{q}}\omega_q(g)D_{1,q}[\Phi;r](g)\right)^2
\end{equation*}
and the \emph{signed asymptotic variance} of the one-level density which is similarly defined by %
\begin{equation*}
\lim_{\substack{q\;\text{prime} \\ q\to+\infty}}%
2%
\smashoperator[r]{%
\sum_{\substack{f\in\prim{\kappa}{q} \\ \epsilon\left(\sym^rf\right)=\epsilon}}%
}%
\omega_q(f)\left(D_{1,q}[\Phi;r](f)-%
2%
\smashoperator[r]{%
\sum_{\substack{g\in\prim{\kappa}{q} \\ \epsilon\left(\sym^rg\right)=\epsilon}}%
}%
\omega_q(g)D_{1,q}[\Phi;r](g)\right)^2
\end{equation*}
when $r$ is odd and $\epsilon=\pm 1$.
\begin{theoint}\label{thm_D}
Let $r\geq 1$ be any integer and $\epsilon=\pm 1$. We assume that hypothesis $\Nice(r,f)$ holds for any prime %
number $q$ and any primitive holomorphic cusp form of level $q$ and even weight $\kappa\geq 2$. If $\nu<1/r^2$ then %
the asymptotic variance of the one-level density is
\begin{equation*}%
2\int_\R\abs{u}\widehat{\Phi}^2(u)\dd u.
\end{equation*}
If $r$ is odd and $\nu<1/(2r(r+2))$ then the signed asymptotic variance of the one-level density is %
\begin{equation*}
2\int_\R\abs{u}\widehat{\Phi}^2(u)\dd u.%
\end{equation*}
\end{theoint}

\subsection{Asymptotic moments of the one-level density}

Last but not least, we compute the \emph{asymptotic $m$-th moment} of the one-level density which is defined by %
\begin{equation*}
\lim_{\substack{q\;\text{prime} \\ q\to+\infty}}%
\smashoperator[r]{%
\sum_{f\in\prim{\kappa}{q}}%
}%
\omega_q(f)\left(D_{1,q}[\Phi;r](f)%
-\sum_{g\in\prim{\kappa}{q}}\omega_q(g)D_{1,q}[\Phi;r](g)\right)^m
\end{equation*}
for any integer $m\geq 1$.%
\begin{theoint}\label{thm_F}
Let $r\geq 1$ be any integer and $\epsilon=\pm 1$. We assume that hypothesis $\Nice(r,f)$ %
holds for any prime number $q$ and any primitive holomorphic cusp form of level $q$ and even weight $\kappa\geq 2$. If $m\nu<4\left/(r(r+2))\right.$ then the asymptotic $m$-th moment of the one-level density is %
\begin{equation*}
\begin{cases}
%\begin{dcases}
0 & \text{if $m$ is odd,}\\%
2\int_\R\abs{u}\widehat{\Phi}^2(u)\dd u\times\frac{m!}{2^{m/2}\left(\frac{m}{2}\right)!} & \text{otherwise.}%
\end{cases}
%\end{dcases}
\end{equation*}
\end{theoint}
\begin{remint}
This result is another evidence for mock-Gaussian behaviour (see \cite{MR2166468,HuRu,HuRu2} for instance).
\end{remint}
\begin{remint}
We compute the first asymptotic moments of the one-level density. These computations allow to compute %
the asymptotic expectation of the first level-densities \cite[\S 1.2]{MR2166468}. We will use the specific case of %
the asymptotic expectation of the two-level density and the asymptotic variance in \S~\ref{sec_twoandvar}.
\end{remint}
Let us sketch the proof of Theorem \ref{thm_F} by explaining the origin of the main term. %
We have to evaluate %
\begin{equation}\label{eq_lasomme}%
\sum_{%
 \substack{%
   0\leq\ell\leq m \\ %
   0\leq\alpha\leq\ell%
          }%
     }%
\binom{m}{\ell}\binom{\ell}{\alpha}R(q)^{\ell-\alpha}%
\Eh[q]\left(P_q^1[\Phi;r]^{m-\ell}P_q^2[\Phi;r]^{\alpha}\right)%
\end{equation}
where $P_q^1[\Phi;r]$ has been defined in \eqref{eq_cneser}, %
\[%
P_q^2[\Phi;r](f)=%
-\frac{2}{\log(q^r)}%
\sum_{j=1}^{r}%
 (-1)^{r-j}%
 \sum_{%
  \substack{%
    p\in\prem\\% 
    p\nmid q%
           }
      }
\lambda_f\left(p^{2j}\right)%
\frac{\log p}{p}\widehat{\Phi}\left(\frac{2\log p}{\log(q^r)}\right)%
\]
and $R(q)$ satisfies %
\[%
R(q)=O\left(\frac{1}{\log q}\right).%
\]
The main term comes from the contribution $\ell=0$ in the sum \eqref{eq_lasomme}. Using a combinatorial lemma, %
we rewrite this main contribution as %
\[%
\frac{(-2)^m}{\log^m{(q^r)}}%
\sum_{s=1}^m%
\sum_{\sigma\in P(m,s)}%
\sum_{\substack{i_1,\dotsc,i_s\\\text{distinct}}}%
\Eh[q]\left(%
\prod_{u=1}^s%
\lambda_f\left(\widehat{p}_{i_{u}}^r\right)^{\varpi^{(\sigma)}_u}%
\right)%
\]
where $P(m,s)$ is the set of surjective functions %
\[%
\sigma \colon \{1,\dotsc,\alpha\}\twoheadrightarrow \{1,\dotsc,s\}%
\]
such that for any $j\in\{1,\dotsc,s\}$, either $\sigma(j)=1$ or there exists $k<j$ such that $\sigma(j)=\sigma(k)+1$ %
and for any $j\in\{1,\dotsc,s\}$ %
\[%
\varpi_j^{(\sigma)}\coloneqq\#\sigma^{-1}(\{j\}).%
\]
$\left(\widehat{p}_i\right)_{i\geq 1}$ stands for the increasing sequence of prime numbers different from $q$. Linearising each $\lambda_f\left(\widehat{p}_{i_{u}}^r\right)^{\varpi^{(\sigma)}_u}$ in terms of %
$\lambda_f\left(\widehat{p}_{i_{u}}^{j_u}\right)$ with $j_u$ runs over integers in $[0,r\varpi^{(\sigma)}_u]$ and using a trace formula to prove that the only $\sigma\in P(m,s)$ leading to a principal contribution satisfy %
$\varpi^{(\sigma)}_j=2$ for any $j\in\{1,\dotsc,s\}$, we have to estimate %
\begin{equation}\label{eq_homo}%
\frac{(-2)^m}{\log^m{(q^r)}}%
\sum_{s=1}^m%
\sum_{\substack{%
\sigma\in P(m,s)\\%
\forall j\in\{1,\dotsc,s\}, \varpi^{(\sigma)}_j=2%
}
}
\sum_{\substack{i_1,\dotsc,i_s\\\text{distinct}}}%
\prod_{u=1}^s%
\frac{%
\log^2{(\widehat{p}_{i_u})}%
}{%%
\widehat{p}_{i_u}%
}%
\widehat{\Phi}^2\left(%
\frac{%
\log\widehat{p}_{i_{u}}%
}{\log{(q^r)}}%
\right).
\end{equation}
This sum vanishes if $m$ is odd since
\[%
\sum_{j=1}^s\varpi_j^{(\sigma)}=m
\]
and it remains to prove the formula for $m$ even. %
In this case, and since we already computed the moment for $m=2$, we deduce from %
\eqref{eq_homo} that the main contribution is %
\[%
\Eh[q](P_q^1[\Phi;r]^2)\times%
\#\left\{%
\sigma\in P(m,m/2) \colon \varpi^{(\sigma)}_j=2\ (\forall j)%
\right\}
\]
and we conclude by computing %
\[%
\#\left\{%
\sigma\in P(m,m/2) \colon \varpi^{(\sigma)}_j=2\ (\forall j)%
\right\}%
=%
\frac{m!}{2^{m/2}\left(\frac{m}{2}\right)!}.
\]
Proving that the other terms lead to error terms is done by implementing similar ideas, but requires -- especially %
for the double products (namely terms implying both $P_q^1$ and $P_q^2$) -- much more combinatorial %
technicalities.

\subsection{Organisation of the paper}
Section \ref{autoproba} contains the automorphic and probabilistic background which is needed to be able to read this paper. In particular, we give here the accurate definition of symmetric power $L$-functions and the properties of Chebyshev polynomials useful in section \ref{momentt}. In section \ref{technical}, we describe the main technical ingredients of this work namely large sieve inequalities for Kloosterman sums and Riemann's explicit formula for symmetric power $L$-functions. In section \ref{one}, some standard facts about symmetry groups are given and the computation of the (signed) asymptotic expectation of the one-level density is done. The computations of the (signed) asymptotic expectation, covariance and variance of the two-level density are done in section \ref{two} whereas the computation of the asymptotic moments of the one-level density is provided in section \ref{momentt}.
Some well-known facts about Kloosterman sums are recalled in appendix \ref{klooster}.

\begin{notations}
We write $\prem$ for the set of prime numbers and the main
parameter in this paper is a prime number $q$, whose name is the level, which
goes to infinity among $\prem$. 
Thus, if $f$ and $g$ are some $\C$-valued functions of the real variable then
the notations $f(q)\ll_{A}g(q)$ or 
$f(q)=O_A(g(q))$ mean that $\abs{f(q)}$ is smaller than a
"constant" which only depends on $A$ times $g(q)$ at least for $q$ a large
enough prime number and similarly, 
$f(q)=o(1)$ means that $f(q)\rightarrow 0$ as $q$ goes to infinity among the
prime numbers. %
We will denote by $\epsilon$ an absolute positive constant whose definition may vary from one line to the next one.
The characteristic function of a set $S$ will be denoted $\un_S$.
\end{notations}

\section{Automorphic and probabilistic background}

\label{autoproba}

\subsection{Automorphic background}\label{sec_autoback}

\subsubsection{Overview of holomorphic cusp forms}

In this section, we recall general facts about holomorphic cusp forms. A reference is \cite{Iw}.
\pa{Generalities} 
We write $\gGamma_{0}(q)$ for the congruence subgroup of level $q$ which acts
on the upper-half plane $\pk$. 
A holomorphic function $f\colon\pk\mapsto\C$ which satisfies
\[%
\forall\begin{pmatrix}
a & b \\
c & d
\end{pmatrix}\in\gGamma_{0}(q),\forall z\in\pk,\quad f\left(\frac{az+b}{cz+d}\right)=(cz+d)^\kappa f(z)
\]
and vanishes at the cusps of $\gGamma_{0}(q)$ is a \emph{holomorphic cusp
  form} of level $q$, even weight $\kappa\geq 2$. We denote by $\cusp{\kappa}{q}$ this space of
  holomorphic cusp forms which is 
equipped with the Peterson inner product
\[%
\scal{f_{1}}{f_{2}}_q\coloneqq%
\int_{\quotientgauche{\gGamma_0(q)}{\pk}}y^{\kappa}f_{1}(z)\overline{f_{2}(z)}\frac{\dd x\dd y}{y^{2}}.
\]
The Fourier expansion at the cusp $\infty$ of any such holomorphic cusp form $f$ is given by
\[%
\forall z\in\pk,\quad f(z)=\sum_{n\geq1}\psi_{f}(n)n^{(\kappa-1)/2}e(nz)
\]
where $e(z)\coloneqq\exp{(2i\pi z)}$ for any complex number $z$. The \emph{Hecke operators} act on $\cusp{\kappa}{q}$ by
\[%%
T_{\ell}(f)(z)
\coloneqq%
\frac{1}{\sqrt{\ell}}\sum_{\substack{ad=\ell\\ (a,q)=1}}\sum_{0\leq b<d}f\left(\frac{az+b}{d}\right)%
\qquad%
\]
for any $z\in\pk$. %
If $f$ is an eigenvector of $T_\ell$, we write $\lambda_f(\ell)$ the corresponding %
eigenvalue. %
We can prove that $T_\ell$ is hermitian if $\ell\geq 1$ is any integer coprime with $q$ and that
\begin{equation}
\label{compo}
T_{\ell_{1}}\circ T_{\ell_{2}}=\sum_{\substack{d\mid(\ell_{1},\ell_{2})\\(d,q)=1}}T_{\ell_{1}\ell_{2}\left/d^{2}\right.}
\end{equation}
for any integers $\ell_1, \ell_2\geq 1$. By Atkin \& Lehner theory \cite{AtLe}, we get a splitting of $\cusp{\kappa}{q}$ 
into $\anc{\kappa}{q}\oplus^{\perp_{\scal{\cdot}{\cdot}_q}}\nouv{\kappa}{q}$ where
\begin{align*}
\anc{\kappa}{q} & \coloneqq   \Vect_{\C}\left\{f(qz), f\in\cusp{\kappa}{1}\right\}\cup\cusp{\kappa}{1}, \\
\nouv{\kappa}{q} & \coloneqq   \left(\anc{\kappa}{q}\right)^{\perp_{\scal{\cdot}{\cdot}_q}}
\end{align*}
where "o" stands for "old" and "n" for "new". Note that $\anc{\kappa}{q}=\{0\}$ if
$\kappa<12$ or $\kappa=14$. 
These two spaces are $T_{\ell}$-invariant for any integer $\ell\geq 1$ coprime with $q$. A
\emph{primitive} cusp form $f\in\nouv{\kappa}{q}$ is an 
eigenfunction of any operator $T_\ell$ for any integer $\ell\geq 1$ coprime with $q$ which is new and arithmetically
normalised 
namely $\psi_{f}(1)=1$. Such an element $f$ is automatically an eigenfunction
of the other Hecke operators 
and satisfies $\psi_{f}(\ell)=\lambda_{f}(\ell)$ for any integer $\ell\geq 1$.
%where $T_\ell(f)=\lambda_f(\ell)f$. 
Moreover, if $p$ is a prime number, define $\alpha_{f}(p)$, $\beta_{f}(p)$ as the complex roots of the quadratic equation
\begin{equation}
\label{quadra}
X^{2}-\lambda_{f}(p)X+\epsilon_{q}(p)=0
\end{equation}
where $\epsilon_{q}$ denotes the trivial Dirichlet character of modulus $q$.
Then it follows from the work of Eichler, Shimura, Igusa and Deligne that
\[%
\abs{\alpha_{f}(p)}, \abs{\beta_{f}(p)}\leq 1 
\]
for any prime number $p$ and so
\begin{equation}
\label{individualh}
\forall\ell\geq 1,\quad 
\abs{\lambda_{f}(\ell)}\leq\tau(\ell).
\end{equation}
The set of primitive cusp forms is denoted by $\prim{\kappa}{q}$. %
It is an orthogonal basis of $\nouv{\kappa}{q}$. %
Let $f$ be a
holomorphic cusp form with Hecke 
eigenvalues $\left(\lambda_{f}(\ell)\right)_{(\ell,q)=1}$. The composition
property \eqref{compo} entails that 
for any integer $\ell_{1}\geq 1$ and for any integer $\ell_{2}\geq 1$ coprime with $q$ the
following multiplicative relations hold:
\begin{align}
\label{compoeigen}
\psi_f(\ell_1)\lambda_f(\ell_2) &= %
\sum_{\substack{d\mid(\ell_1,\ell_2)\\ (d,q)=1}}\psi_f\left(\ell_{1}\ell_{2}\left/d^{2}\right.\right), \\
\psi_f(\ell_1\ell_2) &= %
\sum_{\substack{d\mid(\ell_1,\ell_2)\\ (d,q)=1}}\mu(d)\psi_f\left(\ell_1/d\right)\lambda_f\left(\ell_2/d\right)
\end{align}
and these relations hold for any integers $\ell_{1}, \ell_{2}\geq 1$ if $f$ is primitive. The adjointness relation is
\begin{equation}
\label{adjointness}
\lambda_{f}(\ell)=\overline{\lambda_{f}(\ell)}, \quad \psi_{f}(\ell)=\overline{\psi_{f}(\ell)}
\end{equation}
for any integer $\ell\geq 1$ coprime with $q$ and this remains true for any integer $\ell\geq 1$ if $f$ is primitive.
\pa{Trace formulas} We need two definitions. The harmonic weight associated to any $f$ in $\cusp{\kappa}{q}$ is defined by
\begin{equation}\label{eq_facomeg}%
\omega_q(f)\coloneqq\frac{\fGamma(\kappa-1)}{(4\pi)^{\kappa-1}\scal{f}{f}_q}.%
\end{equation}
For any natural integer $m$ and $n$, the $\Delta_q$-symbol is given by
\begin{equation}\label{eq_deltasymb}
\Delta_q(m,n)\coloneqq\delta_{m,n}+%
2\pi i^{\kappa}\sum_{\substack{c\geq 1 \\ q\mid c}}\frac{S(m,n;c)}{c}J_{\kappa-1}\left(\frac{4\pi\sqrt{mn}}{c}\right)
\end{equation}
where $S(m,n;c)$ is a Kloosterman sum defined in appendix \ref{Kloos} and
$J_{\kappa-1}$ is a Bessel function of first 
kind defined in appendix \ref{Bessel}.The following proposition is \emph{Peterson's trace formula}.
\begin{proposition}\label{prop_orth}
If $\orth{\kappa}{q}$ is any orthogonal basis of $\cusp{\kappa}{q}$ then
\begin{equation}
\label{tr1}
\sum_{f\in\orth{\kappa}{q}}\omega_q(f)\psi_f(m)\psi_f(n)=\Delta_q(m,n)
\end{equation}
for any integers $m$ and $n$.
\end{proposition}
H.~Iwaniec, W.~Luo \& P.~Sarnak proved in \cite{IwLuSa} a useful variation of
Peterson's trace formula which is an average over only primitive cusp
forms. This is more convenient when there are some old forms which is the case
for instance when the weight $\kappa$ is large. Let $\nu$ be the arithmetic function defined by $$\nu(n)\coloneqq n\prod_{p\mid n}\left(1+1/p\right)$$
for any integer $n\geq 1$.
\begin{proposition}[H.~Iwaniec, W.~Luo \& P.~Sarnak~(2001)]%
\label{iwlusatr}%
If $\left(n,q^2\right)\mid q$ and $q\nmid m$ then %
\begin{equation}%
\label{tr2}
\sum_{f\in \prim{\kappa}{q}}\omega_q(f)\lambda_f(m)\lambda_f(n)=%
\Delta_q(m,n)-\frac{1}{q\nu((n,q))}\sum_{\ell\mid q^\infty}\frac{1}{\ell}\Delta_1\left(m\ell^2,n\right).
\end{equation}
\end{proposition}
\begin{remark}
The first term in \eqref{tr2} is exactly the term which appears in \eqref{tr1} %
whereas the second term in \eqref{tr2} %
will be usually very small as an old form comes from a form of level $1$! %
Thus, everything works in practice as if there were no old forms in $\cusp{\kappa}{q}$. %
\end{remark}
\subsubsection{Chebyshev polynomials and Hecke eigenvalues} Let $p\neq q$ a prime number and $f\in\prim{\kappa}{q}$. The multiplicativity %
relation \eqref{compoeigen} leads to %
\[%
\sum_{r\geq 0}\lambda_f(p^r)t^r%
=%
\frac{1}{1-\lambda_f(p)t+t^2}.
\]
It follows that %(\cite[\S 2 et lemme 1]{Ser97})%
\begin{equation}\label{eq_ams}
\lambda_f(p^r)=X_r\left(\lambda_f(p)\right)%
\end{equation}
where the polynomials $X_r$ are defined by their generating series %
\[%
\sum_{r\geq 0}X_r(x)t^r%
=%
\frac{1}{1-xt+t^2}.
\]
They are also defined by
\[%%
X_r(2\cos\theta)=\frac{\sin{((r+1)\theta)}}{\sin{(\theta)}}.
\] 
These polynomials are known as the Chebyshev polynomials of second kind. %
Each $X_r$ has degree $r$, is even if $r$ is even and odd otherwise. %
The family $\{X_r\}_{r\geq 0}$ is a basis for $\Q[X]$, orthonormal %
with respect to the inner product %
\[%
\scalst{P}{Q}\coloneqq\frac{1}{\pi}%
\int_{-2}^{2}P(x)Q(x)\sqrt{1-\frac{x^2}{4}}\dd x.
\] 
In particular, for any integer $\varpi\geq 0$ we have
\begin{equation}\label{eq_lintch}
X_r^\varpi=\sum_{j=0}^{r\varpi}x(\varpi,r,j)X_j
\end{equation}
with
\begin{equation}\label{eq_valx}
x(\varpi,r,j)\coloneqq \scalst{X_r^\varpi}{X_j}%
=\frac{2}{\pi}\int_0^\pi\frac{\sin^\varpi{((r+1)\theta)}\sin{((j+1)\theta)}}{\sin^{\varpi-1}{(\theta)}}\dd\theta.%
\end{equation}
The following relations are useful in this paper %
\begin{equation}
\label{propriox}
x(\varpi,r,j)=\begin{cases}
1 & \text{if $j=0$ and $\varpi$ is even,} \\
0 & \text{if $j$ is odd and $r$ is even,} \\
0 & \text{if $j=0$, $\varpi=1$ and $r\geq 1$.}
\end{cases}
\end{equation}
\subsubsection{Overview of $L$-functions associated to primitive cusp forms} Let $f$ in $\prim{\kappa}{q}$. We define
\[%
L(f,s)\coloneqq\sum_{n\geq 1}\frac{\lambda_f(n)}{n^s}=%
\prod_{p\in\prem}\left(1-\frac{\alpha_f(p)}{p^s}\right)^{-1}\left(1-\frac{\beta_f(p)}{p^s}\right)^{-1}
\]
which is an absolutely convergent and non-vanishing Dirichlet series and Euler product on $\Re{s}>1$ and also
\[%
L_\infty(f,s)\coloneqq\fGamma_{\R}\left(s+(\kappa-1)/2\right)\fGamma_{\R}\left(s+(\kappa+1)/2\right)
\]
where
$\fGamma_{\R}(s)\coloneqq\pi^{-s/2}\fGamma\left(s/2\right)$ as usual.
The function
$$\Lambda(f,s)\coloneqq q^{s/2}L_\infty(f,s)L(f,s)$$
is a
\emph{completed $L$-function} in 
the sense that it satisfies the following \emph{nice} analytic properties:
\begin{itemize}
\item
the function $\Lambda(f,s)$ can be extended to an holomorphic function of order $1$ on $\C$,
\item
the function
$\Lambda(f,s)$ satisfies a functional equation of the shape
\[%
\Lambda(f,s)=i^{\kappa}\epsilon_f(q)\Lambda(f,1-s)
\]
where 
\begin{equation}\label{eq_signe}
\epsilon_f(q)=-\sqrt{q}\lambda_f(q)=\pm 1.
\end{equation}
\end{itemize}

\subsubsection{Overview of symmetric power $L$-functions}\label{sec_sympow}

Let $f$ in $\prim{\kappa}{q}$. For any natural integer $r\geq 1$, the
\emph{symmetric $r$-th power} associated to $f$ is 
given by the following Euler product of degree $r+1$
\[%
L(\sym^rf,s)\coloneqq\prod_{p\in\prem}L_p(\sym^rf,s)
\]
where
\[%
L_p(\sym^rf,s)\coloneqq%
\prod_{i=0}^r\left(1-\frac{\alpha_f(p)^{i}\beta_f(p)^{r-i}}{p^{s}}\right)^{-1}
\]
for any prime number $p$. Let us remark that the local factors of this Euler product may be written as
\[%
L_p(\sym^rf,s)=\prod_{i=0}^r\left(1-\frac{\alpha_f(p)^{2i-r}}{p^s}\right)^{-1}
\]
for any prime number $p\neq q$ and
\[%
L_q(\sym^rf,s)=1-\frac{\lambda_f(q)^r}{q^{s}}=1-\frac{\lambda_f(q^r)}{q^{s}}
\]
as $\alpha_f(p)+\beta_f(p)=\lambda_f(p)$ and
$\alpha_f(p)\beta_f(p)=\epsilon_q(p)$ for any prime number $p$ according to \eqref{quadra}. 
On $\Re{s}>1$, this Euler product is absolutely convergent and
non-vanishing. We also defines \cite[(3.16) and (3.17)]{CoMi} a local 
factor at $\infty$ which is given by a product of $r+1$ Gamma factors namely
\[%
L_\infty(\sym^rf,s)\coloneqq%
\prod_{0\leq a\leq(r-1)/2}%
\fGamma_{\R}\left(s+(2a+1)(\kappa-1)/2\right)%
\fGamma_{\R}\left(s+1+(2a+1)(\kappa-1)/2\right)
\]
if $r$ is odd and
\[%
L_\infty(\sym^rf,s)\coloneqq%
\fGamma_{\R}(s+\mu_{\kappa,r})%
\prod_{1\leq a\leq r/2}%
\fGamma_{\R}\left(s+a(\kappa-1)\right)%
\fGamma_{\R}\left(s+1+a(\kappa-1)\right)
\]
if $r$ is even where
\[%
\mu_{\kappa,r}\coloneqq\begin{cases}
1 & \text{if } r(\kappa-1)/2 \text{ is odd,} \\
0 & \text{otherwise.} \\
\end{cases}
\]
All the local data appearing in these local factors are encapsulated in the following completed $L$-function
\[%
\Lambda(\sym^rf,s)\coloneqq\left(q^{r}\right)^{s/2}L_\infty(\sym^rf,s)L(\sym^rf,s).
\]
Here, $q^r$ is called the arithmetic conductor of $\Lambda(\sym^rf,s)$ and
somehow measures the size of this function. 
We will need some control on the analytic behaviour of this
function. Unfortunately, such information is not currently 
known in all generality. Our main assumption is given in hypothesis $\Nice(r,f)$ page \pageref{hypohypo}. %
Indeed, much more is expected to hold as it is discussed in details in
\cite{CoMi} namely the following assumption 
is strongly believed to be true and lies in the spirit of Langlands program.
\begin{hypCoMi}
There exists an automorphic cuspidal self-dual representation, denoted 
by $\sym^r\pi_f=\otimes'_{p\in\prem\cup\{\infty\}}\sym^r\pi_{f,p}$, of
$GL_{r+1}\left(\A_{\Q}\right)$ whose 
local factors $L\left(\sym^r\pi_{f,p},s\right)$ agree with the local factors $L_p\left(\sym^rf,s\right)$ for any $p$ in $\prem\cup\{\infty\}$.
\end{hypCoMi}
Note that the local factors and the arithmetic conductor in the definition of
$\Lambda\left(\sym^rf,s\right)$ and also the 
sign of its functional equation which all appear without any explanations so
far come from the explicit computations 
which have been done \emph{via} the local Langlands correspondence by J.~Cogdell and
P.~Michel in \cite{CoMi}. 
Obviously, hypothesis $\Nice(r,f)$ is a weak consequence of hypothesis
$\sym^{r}(f)$. For instance, the cuspidality condition in hypothesis $\sym^{r}(f)$ entails the fact that
$\Lambda\left(\sym^rf,s\right)$ is of order $1$ which is crucial 
for us to state a suitable explicit formula. As we will not exploit the power
of automorphic theory in this paper, 
hypothesis $\Nice(r,f)$ is enough for our purpose. In addition, it may happen
that hypothesis $\Nice(r,f)$ is known 
whereas hypothesis $\sym^{r}f$ is not. Let us overview what has been done so
far. 
For any $f$ in $\prim{\kappa}{q}$, hypothesis $\sym^{r}f$ is known for $r=1$
(E.~Hecke), $r=2$ thanks to the work 
of S.~Gelbart and H.~Jacquet \cite{GeJa} and $r=3,4$ from the works of H.~Kim
and F.~Shahidi \cite{KiSh1,KiSh2,Ki}. 

\subsection{Probabilistic background}\label{sec_proba}

The set $\prim{\kappa}{q}$ can be seen as a probability space if
\begin{itemize}
\item
the measurable sets are all its subsets,
\item
the \emph{harmonic probability measure} is defined by
\[%
\muh[q](A)\coloneqq\sumh_{f\in A}1\coloneqq\sum_{f\in A}\omega_q(f)
\]
for any subset $A$ of $\prim{\kappa}{q}$.
\end{itemize}
Indeed, there is a slight abuse here as we only know that
\begin{equation}\label{eq_mesasy}
\lim_{\substack{q\in\prem \\ q\to+\infty}}\muh[q]\left(\prim{\kappa}{q}\right)=1
\end{equation}
(see remark~\ref{rem_moyun}) %
which means that $\muh[q]$ is an ``asymptotic'' probability measure. If $X_q$
is a measurable complex-valued 
function on $\prim{\kappa}{q}$ then it is very natural to compute its \emph{expectation} defined by
\[%
\Eh[q]\left(X_q\right)\coloneqq\sumh_{f\in\prim{\kappa}{q}}X_q(f),
\]
its \emph{variance} defined by
\[%
\Vh[q]\left(X_q\right)\coloneqq\Eh[q]\left(\left(X_q-\Eh[q]\left(X_q\right)\right)^2\right)
\]
and its \emph{$m$-th moments} given by
\[%
\Mh[q,m]\left(X_q\right)\coloneqq\Eh[q]\left(\left(X_q-\Eh[q]\left(X_q\right)\right)^m\right)
\]
for any integer $m\geq 1$. %
If $X\coloneqq\left(X_q\right)_{q\in\prem}$ is a sequence of such measurable
complex-valued functions then we may 
legitimely wonder if the associated complex sequences
\[%
\left(\Eh[q]\left(X_q\right)\right)_{q\in\prem},\quad %
\left(\Vh[q]\left(X_q\right)\right)_{q\in\prem},\quad %
\left(\Mh[q,m]\left(X_q\right)\right)_{q\in\prem}
\]
converge as $q$ goes to infinity among the primes. If yes, the following general notations will be used for their limits
\[%
\Eh[\infty]\left(X\right),\quad %
\Vh[\infty]\left(X\right),\quad%
\Mh[\infty,m]\left(X\right)
\]
for any natural integer $m$. In addition, these potential limits are called
\emph{asymptotic expectation}, 
\emph{asymptotic variance} and \emph{asymptotic $m$-th moments} of $X$ for any
natural integer $m\geq 1$.

For the end of this section, we assume that $r$ is \emph{odd}.
We may remark that the sign of the functional equations of any $L(\sym^rf,s)$ when 
$q$ goes to infinity among the prime numbers and $f$ ranges over $\prim{\kappa}{q}$
is not constant as it depends on $\epsilon_f(q)$.  
Let
\[%
\primeps{\kappa}{q}\coloneqq\left\{f\in \prim{\kappa}{q}, \epsilon(\sym^rf)=\epsilon\right\}
\]
where $\epsilon=\pm 1$. If $f\in\primpair{\kappa}{q}$, then $\sym^rf$ is said
to be \emph{even} whereas it is said to be \emph{odd} if
$f\in\primimpair{\kappa}{q}$. 
It is well-known that
\[%
\lim_{\substack{q\in\prem \\
q\to +\infty}}\muh[q]\left(\left\{f\in\prim{k}{q} \colon \epsilon_f(q)=\epsilon\right\}\right)=\frac{1}{2}.
\]
Since $\epsilon(\sym^rf)$ is $\epsilon_q(f)$ up to a sign depending only on
$\kappa$ and $r$ (by hypothesis $\Nice(r,f)$), it follows that
\begin{equation}\label{eq_sign}
\lim_{\substack{q\in\prem \\
q\to +\infty}}\muh[q]\left(\primeps{\kappa}{q}\right)=\frac{1}{2}.
\end{equation}
For $X_q$ as previous, we can compute its \emph{signed expectation} defined by
\[%
\Eheps[q]\left(X_q\right)\coloneqq%
%\frac{1}{2}\sumh_{f\in\primeps{\kappa}{q}}X_q(f),
2\sumh_{f\in\primeps{\kappa}{q}}X_q(f),
\]
its \emph{signed variance} defined by
\[%
\Vheps[q]\left(X_q\right)\coloneqq\Eheps[q]\left(\left(X_q-\Eheps[q]\left(X_q\right)\right)^2\right)
\]
and its \emph{signed $m$-th moments} given by
\[%
\Mheps[q,m]\left(X_q\right)\coloneqq\Eheps[q]\left(\left(X_q-\Eheps[q]\left(X_q\right)\right)^m\right)
\]
for any natural integer $m\geq 1$. In case of existence, we write
$\Eheps[\infty](X)$, $\Vheps[\infty](X)$ and $\Mheps[\infty,m](X)$ for the
limits which are called \emph{signed asymptotic expectation}, \emph{signed asymptotic variance} and \emph{signed asymptotic moments}. The signed expectation and the expectation are linked through
the formula
\begin{align}
\notag
\Eheps[q](X_q)%
&=%
2\sumh_{f\in\prim{\kappa}{q}}\frac{1+\epsilon\times\epsilon(\sym^rf)}{2}X_q(f)
\\
\label{eq_removeps}
&=%
\Eh[q](X_q)-\epsilon\times\epsilon(\kappa,r)\sqrt{q}\sumh_{f\in\prim{\kappa}{q}}\lambda_f(q)X_q(f).
\end{align}

%..............................
\section{Main technical ingredients of this work}

\label{technical}

\subsection{Large sieve inequalities for Kloosterman sums}

One of the main ingredients in this work is some large sieve inequalities for %
Kloosterman sums which have been %
established by J.-M.~Deshouillers \& H.~Iwaniec in \cite{DeIw} and then %
refined by V.~Blomer, G.~Harcos \& P.~Michel in \cite{BlHaMi}. %
The proof of these large sieve inequalities relies on the %
spectral theory of automorphic forms %
on $GL_2\left(\A_{\Q}\right)$. In particular, the authors have to understand %
the size of the Fourier coefficients of %
these automorphic cusp forms. We have already seen that the size of the %
Fourier coefficients of holomorphic cusp forms %
is well understood \eqref{individualh} but we only have partial results %
on the size of the Fourier %
coefficients of Maass cusp forms which do not come from holomorphic forms. %
We introduce the following hypothesis which measures %
the approximation towards %
the \emph{Ramanujan-Peterson-Selberg conjecture}. %
\begin{ATRPS}\label{hyp_RPS}
If $\pi\coloneqq\otimes'_{p\in\prem\cup\{\infty\}}\pi_{p}$ is any automorphic
cuspidal form on $GL_{2}(\A_\Q)$ with 
local Hecke parameters $\alpha_{\pi}^{(1)}(p)$, $\alpha_{\pi}^{(2)}(p)$ at any
prime number $p$ and 
$\mu_{\pi}^{(1)}(\infty)$, $\mu_{\pi}^{(2)}(\infty)$ at infinity then
\[%
\forall j\in\{1,2\},\quad %
\abs{\alpha_{\pi}^{(j)}(p)}\leq p^{\theta}
\]
for any prime number $p$ for which $\pi_p$ is unramified and
\[%
\forall j\in\{1,2\},\quad%
\abs{%
\Re{\left(\mu_{\pi}^{(j)}(\infty)\right)}
}
\leq\theta
\]
provided $\pi_{\infty}$ is unramified.
\end{ATRPS}
\begin{definition}
We say that\/ $\theta$ is \emph{admissible} if\/ $\Hy_2(\theta)$ is satisfied.
\end{definition}
\begin{remark}
The smallest admissible value of $\theta$ is currently %
$\theta_{0}=\frac{7}{64}$ thanks to the works of H.~Kim, %
F.~Shahidi and %
P.~Sarnak \cite{KiSh2,Ki}. The Ramanujan-Peterson-Selberg conjecture asserts that $0$ is admissible.
\end{remark}
\begin{definition}\label{def_propS}%
Let $T\colon\mathbb{R}^3\to\R^+$ and  $(M,N,C)\in(\mathbb{R}\setminus\{0\})^{3}$, %
we say that a smooth function $h\colon\R^3\to\R^3$ satisfies the property %
$\prp(T;M,N,C)$ if there exists a real number $K>0$ such that %
\begin{multline*}
\forall(i,j,k)\in\N^3, \forall(x_1,x_2,x_3)\in%
\left[\frac{M}{2},2M\right]\times\left[\frac{N}{2},2N\right]\times\left[\frac{C}{2},2C\right], \\
x_1^ix_2^jx_3^k\frac{\partial^{i+j+k} h}%
{\partial x_1^i\partial x_2^j\partial x_3^k}(x_1,x_2,x_3)\leq %
KT(M,N,C)\left(1+\frac{\sqrt{MN}}{C}\right)^{i+j+k}.
\end{multline*}
\end{definition}
With this definition in mind, we are able to write the following proposition %
which is special case of a large sieve inequality adapted from the one %
of Deshouil\-lers \& Iwaniec \cite[Theorem 9]{DeIw} by Blomer, Harcos \& Michel %
\cite[Theorem 4]{BlHaMi}.
\begin{proposition}\label{sieve}%
Let $q$ be some positive integer. %
Let $M, N, C\geq 1$ and $g$ be a smooth function satisfying property $\prp(1;M,N,C)$. %
Consider two sequences of complex numbers $(a_m)_{m\in[M/2,2M]}$ and $(b_n)_{n\in[N/2,2N]}$. %
If \/ $\theta$ is admissible and $MN\ll C^2$ then
\begin{multline}%
\sum_{\substack{c\geq 1 \\ q\mid c}}%
\sum_{m\geq 1}\sum_{n\geq 1}a_mb_n\frac{S(m,\pm n;c)}{c}g(m,n;c)%
\\%
\ll_{\epsilon}%
(qMNC)^\epsilon%
\left(\frac{C^2}{MN}\right)^{\theta}%
\left(%
1+\frac{M}{q}%
\right)^{1/2}%
\left(%
1+\frac{N}{q}%
\right)^{1/2}%
\norm{a}_2\norm{b}_2%
\end{multline}%
for any $\epsilon>0$. %
\end{proposition}%
We shall use a test function.
For any $\nu>0$ let us define $\Schwartz_\nu(\R)$ as the space of even
Schwartz function $\Phi$ %
whose Fourier transform %
\[%
\widehat{\Phi}(\xi)\coloneqq \mathcal{F}[x\mapsto\Phi(x)](\xi)\coloneqq%
\int_{\R}\Phi(x)e(-x\xi)\dd x%
\]
is compactly supported in $[-\nu,+\nu]$. Thanks to the Fourier inversion %
formula: % 
\begin{equation}\label{eq_fouinv}
\Phi(x)=\int_{\R}\widehat{\Phi}(\xi)e(x\xi)\dd x=%
\mathcal{F}[\xi\mapsto\widehat{\Phi}(\xi)](-x),
\end{equation}
such a function $\Phi$ can be extended %
to an entire even function  which satisfies %
\begin{equation}
\label{estim}
\forall s\in\C,\quad %
\Phi(s)\ll_n\frac{\exp{(\nu\abs{\Im{s}})}}{(1+\abs{s})^n}%
\end{equation}
for any integer $n\geq 0$.% 
The version of the large sieve inequality we shall use several times in this paper %
is then the following.%
\begin{corollary}\label{usefulsieve}%
Let $q$ be some prime number, $k_1, k_2>0$ be some integers, $\alpha_1, \alpha_2, \nu$ be some positive real numbers and $\Phi\in\Schwartz_{\nu}(\R)$. Let $h$ be some smooth function satisfying property $\prp(T;M,N,C)$ for any $1\leq M\leq q^{k_1\alpha_1\nu}$, $1\leq N\leq q^{k_2\alpha_2\nu}$ and $C\geq q$. Let $\left(a_{p}\right)_{\substack{p\in\prem\\ p\leq q^{\alpha_1\nu}}}$ and $\left(b_{p}\right)_{\substack{p\in\prem\\ p\leq q^{\alpha_2\nu}}}$ be some complex numbers sequences. If $\;\theta$ is admissible and $\nu\leq2\left/(k_1\alpha_1+k_2\alpha_2)\right.$ then
\begin{multline}%
\sum_{\substack{c\geq 1\\ q\mid c}}%
\sum_{\substack{p_1\in\prem\\ p_1\nmid q}}%
\sum_{\substack{p_2\in\prem\\ p_2\nmid q}}%
a_{p_1}b_{p_2}\frac{S(p_1^{k_1},p_2^{k_2};c)}{c}%
h\left(p_1^{k_1},p_2^{k_2};c\right)%
\widehat{\Phi}\left(\frac{\log p_1}{\log(q^{\alpha_1})}\right)%
\widehat{\Phi}\left(\frac{\log p_2}{\log(q^{\alpha_2})}\right)%
\\%
\ll%
%\\%
q^{\epsilon}%
\sumsh_{\substack{%
1\leq M\leq q^{\nu\alpha_1k_1}\\%
1\leq N\leq q^{\nu\alpha_2k_2}\\%
C\geq q/2%
}}%
\left(1+\sqrt{\frac{M}{q}}\right)%
\left(1+\sqrt{\frac{N}{q}}\right)%
\left(\frac{C^2}{MN}\right)^\theta%
T(M,N,C)%
%\\%
%\times%
\norm{a}_2\norm{b}_2
\end{multline}
where $\sharp$ indicates that the sum is on powers of $\sqrt{2}$. %
The constant implied by the symbol $\ll$ depends at most on $\epsilon$, $k_1$, $k_2$, %
$\alpha_1$, $\alpha_2$ and $\nu$.
\end{corollary}
\begin{proof}%[\proofname{} of corollary \ref{usefulsieve}]
Define $\left(\widehat{a}_m\right)_{m\in\N}$, $\left(\widehat{b}_n\right)_{n\in\N}$ and %
$g(m,n;c)$ by %
\begin{align}%
\widehat{a}_m &\coloneqq a_{m^{1/k_1}}\un_{\prem^{k_1}}(m)\un_{[1,q^{\nu\alpha_1k_1}]}(m)\\%
\widehat{b}_n &\coloneqq b_{n^{1/k_1}}\un_{\prem^{k_1}}(n)\un_{[1,q^{\nu\alpha_1k_1}]}(n)\\%
g(m,n;c) &\coloneqq h(m,n,c)%
\widehat{\Phi}\left(\frac{\log m}{\log(q^{\alpha_1k_1})}\right)%
\widehat{\Phi}\left(\frac{\log n}{\log(q^{\alpha_2k_2})}\right).
\end{align}
Using a smooth partition of unity, as detailed in \S~\ref{unity}, we need to evaluate %
\begin{equation}\label{eq_vtmp}
\sumsh_{\substack{%
1\leq M\leq q^{\nu\alpha_1k_1}\\%
1\leq N\leq q^{\nu\alpha_2k_2}\\%
C\geq q/2%
}}%
T(M,N,C)%
\sum_{\substack{c\geq 1\\ q\mid c}}\sum_{m\geq 1}\sum_{n\geq 1}%
\widehat{a}_m%
\widehat{b}_n%
\frac{S(m,n;c)}{c}%
\frac{g_{M,N,C}(m,n;c)}{T(M,N,C)}.%
\end{equation}%
Since $\nu\leq2\left/(\alpha_1k_1+\alpha_2k_2)\right.$, the first summation is restricted to $MN\ll C^2$ hence, using %
proposition~\ref{sieve}, the quantity in \eqref{eq_vtmp} is %
\begin{equation}%
\ll%
\norm{a}_2\norm{b}_2q^\epsilon%
\sumsh_{\substack{%
1\leq M\leq q^{\nu\alpha_1k_1}\\%
1\leq N\leq q^{\nu\alpha_2k_2}\\%
C\geq q/2%
}}%
T(M,N,C)%
\left(1+\sqrt{\frac{M}{q}}\right)%
\left(1+\sqrt{\frac{N}{q}}\right)%
\left(\frac{C^2}{MN}\right)^\theta.%
\end{equation}
\end{proof}

\subsection{Riemann's explicit formula for symmetric power $L$-functions}\label{sec_explicit}

In this section, we give an analog of Riemann-von Mangoldt's explicit formula %
for symmetric power $L$-functions. Before that, let us %
recall some preliminary facts on % 
zeros of symmetric power $L$-functions which can be found in section 5.3 of %
\cite{IwKo}. Let $r\geq 1$ and $f\in \prim{\kappa}{q}$ for which hypothesis $\Nice(r,f)$ %
holds. All the zeros of $\Lambda(\sym^rf,s)$ are %
in the critical strip $\{s\in\C\colon 0<\Re{s}<1\}$. The multiset of the zeros %
of $\Lambda(\sym^rf,s)$ counted with multiplicities is given by %
\[%
\left\{%
\rho_{f,r}^{(j)}=\beta_{f,r}^{(j)}+i\gamma_{f,r}^{(j)} \colon j\in\mathcal{E}(f,r)
\right\}
\]
where
\[%
\mathcal{E}(f,r)\coloneqq%
\begin{cases}
\Z & \text{if $\sym^rf$ is odd}\\
\Z\setminus\{0\} & \text{if $\sym^rf$ is even.}
\end{cases}
\]
and
\begin{align*}
\beta_{f,r}^{(j)} & =  \Re{\rho_{f,r}^{(j)}}, \\
\gamma_{f,r}^{(j)} & =  \Im{\rho_{f,r}^{(j)}}
\end{align*}
for any $j\in\mathcal{E}(f,r)$.
We enumerate the zeros such that
\begin{enumerate}
\item the sequence $j\mapsto\gamma_{f,r}^{(j)}$ is increasing
\item we have $j\geq 0$ if and only if  $\gamma_{f,r}^{(j)}\geq 0$
\item we have $\rho_{f,r}^{(-j)}=1-\rho_{f,r}^{(j)}$.
\end{enumerate}
Note that if $\rho_{f,r}^{(j)}$ is a
zero of $\Lambda(\sym^rf,s)$ then 
$\overline{\rho_{f,r}^{(j)}}$, $1-\rho_{f,r}^{(j)}$ and $1-\overline{\rho_{f,r}^{(j)}}$ are also
some zeros of $\Lambda(\sym^rf,s)$. In addition, remember that if $\sym^rf$ is odd then the functional equation of $L(\sym^rf,s)$ evaluated at the critical point $s=1/2$ provides a trivial zero denoted by $\rho_{f,r}^{(0)}$. It can be shown \cite[Theorem 5.8]{IwKo} 
that the number of zeros $\Lambda(\sym^rf,s)$ satisfying $\abs{\gamma_{f,r}^{(j)}}\leq T$ is
\begin{equation}\label{eq_nbzero}
\frac{T}{\pi}\log{\left(\frac{q^rT^{r+1}}{(2\pi e)^{r+1}}\right)}+O\left(\log(qT)\right)
\end{equation}
as $T\geq 1$ goes to infinity. We state now the \emph{Generalised Riemann
  Hypothesis} which is the main conjecture 
about the horizontal distribution of the zeros of $\Lambda(\sym^rf,s)$ in the critical strip.
\begin{hGRH}
For any prime number $q$ and any $f$ in $\prim{\kappa}{q}$, all the zeros of
$\Lambda(\sym^rf,s)$ lie on the critical 
line $\left\{s\in\C \colon \Re{s}=1/2\right\}$ namely $\beta_{r,f}^{(j)}=1/2$ for any $j\in\mathcal{E}(f,r)$.
\end{hGRH}
\begin{remark}
We \emph{do not} use this hypothesis in our proofs.
\end{remark}
Under hypothesis $\GRH(r)$, it can be shown that the number of zeros of %
the function %
$\Lambda(\sym^rf,s)$ satisfying % 
$\abs{\gamma_{f,r}^{(j)}}\leq 1$ is given by %
\[%
\frac{1}{\pi}\log{\left(q^r\right)}(1+o(1))
\]
as $q$ goes to infinity. Thus, the spacing between two consecutive zeros with imaginary part in $[0,1]$ is roughly of size
\begin{equation}
\label{meanspacing}
\frac{2\pi}{\log{\left(q^r\right)}}.
\end{equation} 
We aim at studying the local distribution of the zeros of $\Lambda(\sym^rf,s)$
in a neighborhood of the real axis of size $1/\log q^r$ since in such a neighborhood, we expect
to catch only few zeros (but without being able to say that we catch only one\footnote{We refer to 
Miller \cite{mil02a} and Omar \cite{oma00} for works related to the ``first'' zero.}).
Hence, we normalise the zeros by defining
\[%
\widehat{\rho}_{f,r}^{(j)}%
\coloneqq %
\frac{\log{\left(q^r\right)}}{2i\pi}\left(%
\beta_{f,r}^{(j)}-\frac{1}{2}+i\gamma_{f,r}^{(j)}
\right).
\]
Note that 
\[%
\widehat{\rho}_{f,r}^{(-j)}=-\widehat{\rho}_{f,r}^{(j)}.
\]
\begin{definition}
Let $f\in\prim{\kappa}{q}$ for which hypothesis\/ $\Nice(r,f)$ %
holds and let $\Phi\in\Schwartz_\nu(\R)$. The \emph{one-level density} %
(relatively to $\Phi$) of\/ $\sym^rf$ is %
\begin{equation}
\label{eq_defdens}
D_{1,q}[\Phi;r](f)%
\coloneqq%
\sum_{j\in\mathcal{E}(f,r)}%
\Phi\left(\widehat{\rho}_{f,r}^{(j)}\right).
\end{equation}
\end{definition}

To study $D_{1,q}[\Phi;r](f)$ for any $\Phi\in\Schwartz_\nu(\R)$, we transform this sum over zeros into a sum over primes in the next
proposition. In other words, we establish an explicit formula for symmetric power $L$-functions. Since the proof is classical, we refer to
\cite[\S 4]{IwLuSa} or \cite[\S 2.2]{Gu} which present a method that has just to be adapted to our setting.
\begin{proposition}
\label{explicit}
Let $r\geq 1$ and $f\in \prim{\kappa}{q}$ for which hypothesis $\Nice(r,f)$ holds and let $\Phi\in\Schwartz_\nu(\R)$. We have
\[%
D_{1,q}[\Phi;r](f)=%
E[\Phi;r]+P_q^1[\Phi;r](f)+\sum_{m=0}^{r-1}(-1)^mP_q^2[\Phi;r,m](f)+O\left(\frac{1}{\log{\left(q^r\right)}}\right)
\]
where
\begin{align*}
E[\Phi;r] & \coloneqq  \widehat{\Phi}(0)+\frac{(-1)^{r+1}}{2}\Phi(0), \\
P_{q}^1[\Phi;r](f) & \coloneqq %
-\frac{2}{\log{\left(q^r\right)}}\sum_{\substack{p\in\prem \\ p\nmid q}}%
\lambda_{f}\left(p^r\right)\frac{\log{p}}{\sqrt{p}}\widehat{\Phi}\left(\frac{\log{p}}{\log{\left(q^r\right)}}\right), \\
P_q^2[\Phi;r,m](f) & \coloneqq  %
-\frac{2}{\log{\left(q^r\right)}}\sum_{\substack{p\in\prem \\ p\nmid q}}%
\lambda_f\left(p^{2(r-m)}\right)\frac{\log{p}}{p}\widehat{\Phi}\left(\frac{2\log{p}}{\log{\left(q^r\right)}}\right)
\end{align*}
for any integer $m\in\{0,\ldots,r-1\}$.
\end{proposition}
%..........................................
\subsection{Contribution of the old forms}
In this short section, we prove the following useful lemmas.
%..............................................
\begin{lemma}\label{lem_delun}
Let $p_1$ and $p_2\neq q$ be some prime numbers and $a_1$, $a_2$, $a$ be some nonnegative integers. Then %
\[%%
\sum_{\ell\mid q^{\infty}}\frac{\Delta_1\n(\ell^2p_1^{a_1},p_2^{a_2}q^a)}{\ell}%
\ll%
\frac{1}{q^{a/2}}%
\]
the implied constant depending only on $a_1$ and $a_2$.
\end{lemma}
\begin{proof}
Using proposition~\ref{prop_orth} and the fact that $\orth{\kappa}{1}=\prim{\kappa}{1}$, we write %
\begin{align}
\Delta_1\n(\ell^2p_1^{a_1},p_2^{a_2}q^a)%
&=%
\sumh_{f\in\prim{\kappa}{1}}\lambda_f\n(\ell^2p_1^{a_1})\lambda_f\n(p_2^{a_2}q^a)\\
&\ll%
\label{eq_ici}
\sumh_{f\in\prim{\kappa}{1}}\abs{\lambda_f\n(\ell^2p_1^{a_1})}\cdot\abs{\lambda_f\n(p_2^{a_2})}\cdot\abs{\lambda_f(q^a)}.
\end{align}
By Deligne's bound~\eqref{individualh} we have
\begin{equation}\label{eq_sansq}
\abs{\lambda_f\n(\ell^2p_1^{a_1})}\cdot\abs{\lambda_f\n(p_2^{a_2})}%
\leq%
\tau(\ell^2p_1^{a_1})\tau(p_2^{a_2})%
\leq (a_1+1)(a_2+2)\tau(\ell^2).
\end{equation}
By the multiplicativity relation~\eqref{compoeigen} and the value of the sign of the functional equation~\eqref{eq_signe}, we have
\begin{equation}\label{eq_avecq}
\abs{\lambda_f(q^a)}\ll%
\frac{1}{q^{a/2}}.%
\end{equation}
We obtain the result by reporting \eqref{eq_avecq} and \eqref{eq_sansq} in \eqref{eq_ici} and by using \eqref{eq_mesasy} and
\[%
\sum_{\ell\mid q^{\infty}}\frac{\tau(\ell^2)}{\ell}=\frac{1+1/q}{(1-1/q)^2}\ll 1.
\]
\end{proof}
\begin{lemma}\label{deltaestimate}
Let $m,n\geq 1$ be some coprime integers. Then,
\[%
\Delta_q(m,n)-\delta(m,n)\ll%
%\begin{dcases}
\begin{cases}
\frac{(mn)^{1/4}}{q}\log\left(\frac{mn}{q^2}\right) & \text{if $mn>q^2$}\\
\frac{(mn)^{( \kappa-1)/2}}{q^{\kappa-1/2}} %
\leq \frac{(mn)^{1/4}}{q} & \text{if $mn\leq q^2$.}
%\end{dcases}
\end{cases}
\]
\end{lemma}
\begin{proof}
This is a direct consequence of the Weil-Estermann bound \eqref{weil} %
and lem\-ma \ref{lem_picard}.
\end{proof}
\begin{corollary}\label{lem_sumlambdafq}
For any prime number $q$, we have
\[%
\sqrt{q}\sumh_{f\in\prim{\kappa}{q}}\lambda_f(q)\ll\frac{1}{q^{\delta_\kappa}}
\]
where
\[%%
\delta_\kappa\coloneqq%
%\begin{dcases*}
\begin{cases}
\frac{\kappa-1}{2} & \text{if $\kappa\leq 10$ or $\kappa=14$} \\ %
\frac{5}{2} & \text{otherwise.}%
%\end{dcases*}
\end{cases}
\]
\end{corollary}
\begin{proof}[\proofname{} of corollary~\ref{lem_sumlambdafq}]
Let $\mathcal{K}=\{\kappa\in 2\N \colon 2\leq\kappa\leq 14,\, \kappa\neq 12\}$. %
By proposition~\ref{iwlusatr}, we have
\begin{equation}\label{eq_termzero}%
\sumh_{f\in\prim{\kappa}{q}}\lambda_f(q)=%
\Delta_q(1,q)-%
\frac{\delta(\kappa%
\notin%
\mathcal{K})}{q\nu(q)}%
\sum_{\ell\mid q^\infty}\frac{\Delta_1(\ell^2,q)}{\ell}.
\end{equation}
The term $\delta(\kappa\notin\mathcal{K})$ comes from
proposition~\ref{prop_orth} with the fact that %
there is no cusp forms of %
weight $\kappa\in\mathcal{K}$ and level $1$. Lemma~\ref{deltaestimate} gives %
\begin{equation}\label{eq_termun}
\Delta_q(1,q)\ll \frac{1}{q^{\kappa/2}}
\end{equation}
and lemma~\ref{lem_delun} gives %
\begin{equation}\label{eq_termdeux}%
\sum_{\ell\mid q^\infty}\frac{\Delta_1(\ell^2,q)}{\ell}\ll\frac{1}{%
\sqrt{q}%
}.
\end{equation}
Since $\nu(q)>q$, the result follows from reporting \eqref{eq_termun} and \eqref{eq_termdeux} in \eqref{eq_termzero}. 
\end{proof}
\begin{remark}\label{rem_moyun}
In a very similar fashion, one can prove that
\begin{equation}\label{eq_moyun}
\muh[q]\left(\prim{\kappa}{q}\right)%
=\Eh[q](1)%
=1+O\left(\frac{1}{q^{\gamma_\kappa}}\right).
\end{equation}
where
\[%
\gamma_\kappa\coloneqq%
%\begin{dcases*}
\begin{cases}
\kappa-\frac{1}{2} & \text{if $\kappa\leq 10$ or $\kappa=14$} \\ %
1 & \text{otherwise.}%
%\end{dcases*}
\end{cases}
\]
Corollary \ref{lem_sumlambdafq}, \eqref{eq_moyun} and \eqref{eq_removeps} %
imply %
\begin{equation}\label{eq_mbun}
\Eheps[q](1)%
=1+O\left(\frac{1}{q^{\beta_\kappa}}\right)
\end{equation}
where %
\[%
\beta_\kappa\coloneqq
\begin{cases}
%\begin{dcases*}
\frac{\kappa-1}{2} & \text{if $\kappa\leq 10$ or $\kappa=14$} \\ %
1 & \text{otherwise.}%
%\end{dcases*}
\end{cases}
\]
\end{remark}
%..............................................
A direct consequence of lemma~\ref{lem_delun} is the following one.
\begin{lemma}\label{lem_old}
Let $\alpha_1,\alpha_2,\beta_1,\beta_2,\gamma_1,\gamma_2,w$ be some nonnegative real numbers.
Let $\Phi_1$ and $\Phi_2$ be in $\Schwartz_\nu(\R)$.
Then,
\begin{multline*}
\sum_{\substack{p_1\in\prem\\ p_1\nmid q}}%
\sum_{\substack{p_2\in\prem\\ p_2\nmid q}}%
\frac{\log p_1}{p_1^{\alpha_1}}\frac{\log p_2}{p_2^{\alpha_2}}%
\widehat{\Phi_1}\left(\frac{\log p_1}{\log{\left(q^{\beta_1}\right)}}\right)%
\widehat{\Phi_2}\left(\frac{\log p_2}{\log{\left(q^{\beta_2}\right)}}\right)%
\sum_{\ell\mid q^\infty}\frac{\Delta_1(\ell^2p_1^{\gamma_1},p_2^{\gamma_2}q^w)}{\ell}%
\\
\ll%
q^{\delta\nu-w/2+\varepsilon}
\end{multline*}
with $\delta$ given in table~\ref{tab_vald}.
\begin{table}[H]
%\begin{center}
\setlength{\extrarowheight}{4pt}
\begin{tabular}{|c|c|c|}
\hline
\backslashbox{$\alpha_2$}{$\alpha_1$} & $]0,1]$ & $[1,+\infty[$\\ 
\hline
$]0,1]$ & $\beta_1(1-\alpha_1)+\beta_2(1-\alpha_2)$ & $\beta_2(1-\alpha_2)$ \\
\hline
$[1,+\infty[$ & $\beta_1(1-\alpha_1)$ & $0$ \\
\hline
\end{tabular}
\caption{Values of $\delta$}\label{tab_vald}
%\end{center}
\end{table}
\end{lemma}
%..........................................
\section{Linear statistics for low-lying zeros}
\label{one}
\subsection{Density results for families of $L$-functions}\label{sec_densres}

We briefly recall some well-known features that can be found in
\cite{IwLuSa}. 
Let $\mathcal{F}$ be a family of $L$-functions indexed by the arithmetic conductor namely
\[%
\mathcal{F}=\bigcup_{Q\geq 1}\mathcal{F}(Q)
\]
where the arithmetic conductor of any $L$-function in $\mathcal{F}(Q)$ is of
order $Q$ in the logarithmic scale. 
It is expected that there is a symmetry group $G(\mathcal{F})$ of matrices of
large rank endowed with a probability 
measure which can be associated to $\mathcal{F}$ such that the low-lying zeros
of the $L$-functions in $\mathcal{F}$ 
namely the non-trivial zeros of height less than $1/\log{Q}$ are
distributed like the eigenvalues of the 
matrices in $G(\mathcal{F})$. In other words, there should exist a symmetry
group $G(\mathcal{F})$ such that for 
any $\nu>0$ and any $\Phi\in\Schwartz_\nu(\R)$,
\begin{multline*}
\lim_{Q\to+\infty}\frac{1}{\mathcal{F}(Q)}%
\sum_{\pi\in\mathcal{F}(Q)}\sum_{\substack{%
0\leq\beta_\pi\leq 1 \\
\gamma_\pi\in\R \\
L\left(\pi,\beta_\pi+i\gamma_\pi\right)=0%
}}
\Phi\left(\frac{\log{Q}}{2i\pi}\left(\beta_\pi-\frac{1}{2}+i\gamma_\pi\right)\right) \\
=%
\int_{\R}\Phi(x)W_1(G(\mathcal{F}))(x)\dd x
\end{multline*}
where $W_1(G(\mathcal{F}))$ is the one-level density of the eigenvalues of
$G(\mathcal{F})$. In this case, 
$\mathcal{F}$ is said to be of \emph{symmetry type} $G(\mathcal{F})$ and we
said that we proved a \emph{density result} 
for $\mathcal{F}$. For instance, the following densities are determined in \cite{KaSa}:
\begin{align*}
W_1(SO(\mathrm{even}))(x) &= 1+\frac{\sin{(2\pi x)}}{2\pi x},\\
W_1(O)(x) &= 1+\frac{1}{2}\delta_0(x),\\
W_1(SO(\mathrm{odd}))(x) &= 1-\frac{\sin{(2\pi x)}}{2\pi x}+\delta_0(x),\\
W_1(Sp)(x) &= 1-\frac{\sin{(2\pi x)}}{2\pi x}
\end{align*}
where $\delta_0$ is the Dirac distribution at $0$. According to Plancherel's formula,
\[%
\int_{\R}\Phi(x)W_1(G(\mathcal{F}))(x)\dd x=\int_{\R}\widehat{\Phi}(x)\widehat{W}_1(G(\mathcal{F}))(x)\dd x
\]
and we can check that
\begin{align*}
\widehat{W}_1(SO(\mathrm{even}))(x) &= \delta_0(x)+\frac{1}{2}\eta(x),\\
\widehat{W}_1(O)(x) &= \delta_0(x)+\frac{1}{2},\\
\widehat{W}_1(SO(\mathrm{odd}))(x) &= \delta_0(x)-\frac{1}{2}\eta(x)+1,\\
\widehat{W}_1(Sp)(x) &= \delta_0(x)-\frac{1}{2}\eta(x)
\end{align*}
where
\[%
\eta(x)\coloneqq
\begin{cases}
1 & \text{if $\abs{x}<1$,} \\
\frac{1}{2} & \text{if $x=\pm 1$,} \\
0 & \text{otherwise.}
\end{cases}
\]
As a consequence, if we can only prove a density result for $\nu\leq 1$, the
three orthogonal densities are 
indistinguishable although they are distinguishable from $Sp$. Thus, the challenge is to pass the natural barrier $\nu=1$.

\subsection{Asymptotic expectation of the one-level density}

The aim of this part is to prove a density result for the family
\[%
\mathcal{F}_r\coloneqq\bigcup_{q\in\prem}\left\{L(\sym^rf,s), f\in \prim{\kappa}{q}\right\}
\]
for any $r\geq 1$ which consists in proving the existence and computing the
asymptotic 
expectation $\Eh[\infty]\left(D_{1}[\Phi;r]\right)$ of
$D_1[\Phi;r]\coloneqq\left(D_{1,q}[\Phi;r]\right)_{q\in\prem}$ 
for any $r\geq 1$ and for $\Phi$ in $\Schwartz_\nu(\R)$ with $\nu>0$ as large
as possible in order to be able to 
distinguish between the three orthogonal densities if $r$ is small enough.
Recall that $E[\Phi;r]$ has been defined in proposition \ref{explicit}.
\begin{theorem}
\label{density1}
Let $r\geq 1$ and $\Phi\in\Schwartz_\nu(\R)$. We assume that hypothesis
$\Nice(r,f)$ holds for any prime number $q$ 
and any $f\in \prim{\kappa}{q}$ and also that $\theta$ is admissible. Let
\[%
\nu_{1,\mathrm{max}}(r,\kappa,\theta)\coloneqq\left(1-\frac{1}{2(\kappa-2\theta)}\right)\frac{2}{r^2}.
\]
If $\nu<\nu_{1,\mathrm{max}}(r,\kappa,\theta)$ then
\[%
\Eh[\infty]\left(D_{1}[\Phi;r]\right)=E[\Phi;r].
\]
\end{theorem}
\begin{remark}
We remark that
\begin{alignat}{2}
\label{eq_pqcn}
\nu_{1,\mathrm{max}}(r,\kappa,\theta_0) &= \left(1-\frac{16}{32\kappa-7}\right)\frac{2}{r^2}& \geq \frac{82}{57r^2}, \\
\nu_{1,\mathrm{max}}(r,\kappa,0) &= \left(1-\frac{1}{2\kappa}\right)\frac{2}{r^2}& \geq \frac{3}{2r^2} 
\end{alignat}
and thus $\nu_{1,\mathrm{max}}(1,\kappa,\theta_0)>1$ whereas $\nu_{1,\mathrm{max}}(r,\kappa,\theta_0)\leq 1$ for any $r\geq 2$.
\end{remark}
\begin{remark}\label{rem_symtyp}
Note that
\[%
E[\Phi;r]=\int_{\R}\widehat{\Phi}(x)\left(\delta_0(x)+\frac{(-1)^{r+1}}{2}\right)\dd x.
\]
Thus, this theorem reveals that the symmetry type of $\mathcal{F}_r$ is
\[%
G(\mathcal{F}_r)=
\begin{cases}
Sp & \text{if $r$ is even,} \\
O & \text{if $r=1$,} \\
SO(\mathrm{even}) \text{ or } O \text{ or } SO(\mathrm{odd}) & \text{if $r\geq 3$ is odd.}
\end{cases}
\]
Some additional comments are given in remark \ref{remark4} page \pageref{remark4}.
\end{remark}
\begin{proof}[\proofname{} of theorem \ref{density1}.] The proof is detailed and will
  be a model for the next density results. 
According to proposition \ref{explicit} %
and \eqref{eq_moyun}
, we have
\begin{multline}
\label{explicitaverage}
\Eh[q]\left(D_{1,q}[\Phi;r]\right)=%
E[\Phi;r]%
+\Eh[q]\left(P_q^1[\Phi;r]\right) \\
+\sum_{m=0}^{r-1}(-1)^m\Eh[q]\left(P_q^2[\Phi;r,m]\right)+O\left(\frac{1}{\log{\left(q^r\right)}}\right).
\end{multline}
The first term in \eqref{explicitaverage} is the main term given in the
theorem. We now estimate the second 
term of \eqref{explicitaverage} \emph{via} the trace formula given in proposition \ref{iwlusatr}.
\begin{equation}\label{eq_pu}
\Eh[q]\left(P_q^1[\Phi;r]\right)=%
\PP{q,\mathrm{new}}1[\Phi;r]+\PP{q,\mathrm{old}}1[\Phi;r]
\end{equation}
where
\begin{align*}
\PP{q,\mathrm{new}}1[\Phi;r] &= %
-\frac{2}{\log{\left(q^r\right)}}\sum_{\substack{p\in\prem \\ p\nmid q}}\Delta_q(p^r,1)\frac{\log{p}}{\sqrt{p}}%
\widehat{\Phi}\left(\frac{\log{p}}{\log{\left(q^r\right)}}\right), \\
\PP{q,\mathrm{old}}1[\Phi;r] &= \frac{2}{q\log{\left(q^r\right)}}%
\sum_{\ell\mid q^\infty}\frac{1}{\ell}\sum_{\substack{p\in\prem \\ p\nmid q}}\Delta_1(p^r\ell^2,1)%
\frac{\log{p}}{\sqrt{p}}\widehat{\Phi}\left(\frac{\log{p}}{\log{\left(q^r\right)}}\right).
\end{align*}
Let us estimate the new part which can be written as
\begin{multline*}
\PP{q,\mathrm{new}}1[\Phi;r]=%
-\frac{2%
(2\pi i^\kappa)
}{\log{\left(q^r\right)}}\sum_{\substack{c\geq 1 \\ q\mid c}}%
\sum_{p\in\prem}\left(\frac{\log{p}}{\sqrt{p}}\delta_{q\nmid p}\un_{\left[1,q^{r\nu}\right]}(p)\right)%
\frac{S(p^r,1;c)}{c} \\
\times J_{\kappa-1}\left(\frac{4\pi\sqrt{p^r}}{c}\right)\widehat{\Phi}\left(\frac{\log{p}}{\log{\left(q^{r}\right)}}\right).
\end{multline*}
Thanks to \eqref{bessel}, the function %
\[%
h(m;c)\coloneqq J_{\kappa-1}\left(\frac{4\pi\sqrt{m}}{c}\right)%
\]%
satisfies hypothesis $\prp(T;M,1,C)$ with %
\[%
T(M,1,C)=\left(1+\frac{\sqrt{M}}{C}\right)^{1/2-\kappa}%
\left(\frac{\sqrt{M}}{C}\right)^{\kappa-1}.
\]%
Hence, if $\nu\leq 2/r^2$ then corollary \ref{usefulsieve} leads to %
\begin{align}%
\label{eq_plusgrand}%
\PP{q,\mathrm{new}}1[\Phi;r]%
&\ll_\epsilon%
q^\epsilon%
\sumsh_{\substack{1\leq M\leq q^{\nu r^2}\\ C\geq q/2}}%
\left(1+\sqrt{\frac{M}{q}}\right)%
\left(\frac{\sqrt{M}}{C}\right)^{\kappa-1-2\theta}%
\\%
&\ll_\epsilon%
q^\epsilon%
\sumsh_{1\leq M\leq q^{\nu r^2}}%
\left(\frac{M^{\frac{\kappa-1}{2}-\theta}}{q^{\kappa-1-2\theta}}+%
\frac{M^{\frac{\kappa}{2}-\theta}}{q^{\kappa-\frac{1}{2}-2\theta}}\right)%
\end{align}%
thanks to \eqref{dyadic2}.
Summing over $M$ \emph{via} \eqref{dyadic1} leads to
\begin{equation}\label{eq_pun}
\PP{q,\mathrm{new}}1[\Phi;r]\ll_\epsilon %
q^{\left(\frac{\kappa-1}{2}-\theta\right)(r^2\nu-2)+\epsilon}+%
q^{\left(\frac{\kappa}{2}-\theta\right)r^2\nu-\left(\kappa-\frac{1}{2}-2\theta\right)+\epsilon}
\end{equation}
which is an admissible error term if $\nu<\nu_{1,\mathrm{max}}(r,\kappa,\theta)$. %
According to lemma \ref{lem_old} (with $\alpha_2=+\infty$) we have %
\begin{equation}\label{eq_puo}
\PP{q,\mathrm{old}}1[\Phi;r]\ll_\epsilon q^{\frac{r\nu}{2}-1+\epsilon}
\end{equation}
which is an admissible error term if $\nu<2/r$. Reporting
\eqref{eq_pun} and \eqref{eq_puo} in \eqref{eq_pu} we obtain
\begin{equation}\label{eq_psum}
\Eh[q]\left(P^1_q[\Phi;r]\right)\ll\frac{1}{q^{\delta_1}}
\end{equation}
for some $\delta_1>0$ (depending on $\nu$ and $r$) as soon as $\nu<\nu_{1,\mathrm{max}}(r,\kappa,\theta)$.
We now estimate the
third term of \eqref{explicitaverage}. 
If $0\leq m\leq r-1$ then the trace formula given in proposition \ref{iwlusatr} implies that
\begin{equation}\label{eq_pd}
\Eh[q]\left(P_q^2[\Phi;r,m]\right)=\PP{q,\mathrm{new}}2[\Phi;r,m]+\PP{q,\mathrm{old}}2[\Phi;r,m]
\end{equation}
where
\begin{align*}
\PP{q,\mathrm{new}}2[\Phi;r,m] &= %
-\frac{2}{\log{\left(q^r\right)}}\sum_{\substack{p\in\prem \\ p\nmid q}}%
\Delta_q\left(p^{2(r-m)},1\right)\frac{\log{p}}{p}\widehat{\Phi}%
\left(\frac{\log{\left(p^2\right)}}{\log{\left(q^r\right)}}\right), \\
\PP{q,\mathrm{old}}2[\Phi;r,m] &= \frac{2}{q\log{\left(q^r\right)}}%
\sum_{\ell\mid q^\infty}\frac{1}{\ell}\sum_{\substack{p\in\prem \\ p\nmid q}}%
\Delta_1\left(p^{2(r-m)}\ell^2,1\right)\frac{\log{p}}{p}\widehat{\Phi}%
\left(\frac{\log{\left(p^2\right)}}{\log{\left(q^r\right)}}\right).
\end{align*}
Let us estimate the new part which can be written as
\begin{multline*}
\PP{q,\mathrm{new}}2[\Phi;r,m]=-\frac{2%
(2\pi i^\kappa)%
}{\log{\left(q^r\right)}}\sum_{\substack{c\geq 1 \\
q\mid c}}\sum_{p\in\prem}\left(\frac{\log{p}}{\sqrt{p}}\delta_{q\nmid p}\un_{\left[1,q^{\frac{r\nu}{2}}\right]}(p)\right)\frac{S\left(p^{2(r-m)},1;c\right)}{c} \\
\times \frac{1}{\sqrt{p}}J_{\kappa-1}\left(\frac{4\pi\sqrt{p^{2(r-m)}}}{c}\right)%
\widehat{\Phi}\left(%
\frac{\log p}{\log q^{r/2}}
\right).
\end{multline*}
The function %
\[%
h(m,c)\coloneqq J_{\kappa-1}\left(\frac{4\pi\sqrt{m}}{c}\right)\times\frac{1}{m^{1/(4(r-m))}}%
\]%
satisfies hypothesis $\prp(T;M,1,C)$ with %
\[%
T(M,1,C)=\left(1+\frac{\sqrt{M}}{C}\right)^{1/2-\kappa}%
\left(\frac{\sqrt{M}}{C}\right)^{\kappa-1}%
\frac{1}{M^{1/(4(r-m))}}.
\]%
Hence, if $\nu\leq 2/r^2$ then corollary~\ref{usefulsieve} leads to %
\[%
\PP{q,\mathrm{new}}2[\Phi;r,m]%
\ll_\epsilon%
q^\epsilon%
\sumsh_{\substack{M\leq q^{\nu r(r-m)}\\ C\geq q/2}}%
\frac{1}{(M)^{1/(4r-4m)}}%
\left(\frac{\sqrt{M}}{C}\right)^{\kappa-1-2\theta}%
\left(1+\sqrt{\frac{M}{q}}\right).%
\]
This is smaller than the bound given in %
\eqref{eq_plusgrand} and hence %
is an admissible error term if $\nu<\nu_{1,\mathrm{max}}(r,\kappa,\theta)$. % 
 According to lemma \ref{lem_old}, we have %
 \begin{equation}\label{eq_pdo}
 \PP{q,\mathrm{old}}2[\Phi;r]\ll_\epsilon %
 q%
^{-1+\epsilon}.
\end{equation}
We obtain
\begin{equation}\label{eq_psdm}
\Eh[q]\left(P^2_q[\Phi;r%
,m%
]\right)\ll\frac{1}{q^{\delta_2}}
\end{equation}
for some $\delta_2>0$ (depending on $\nu$ and $r$) as soon as $\nu<\nu_{1,\mathrm{max}}(r,\kappa,\theta)$.
Finally, reporting \eqref{eq_psdm} and \eqref{eq_psum} in %
\eqref{explicitaverage}, we get % 
\begin{equation}\label{eq_mhe}
\Eh[q]\left(D_{1,q}[\Phi;r]\right)=%
E[\Phi;r]+O\left(\frac{1}{\log q}\right).
\end{equation}
\end{proof}

%........
\subsection{Signed asymptotic expectation of the one-level density}
In this part, we prove some density results for
subfamilies of $\mathcal{F}_r$ on which the sign of the functional equation remains constant.
The two subfamilies are defined by
\[%%
\mathcal{F}_r^{\epsilon}\coloneqq\bigcup_{q\in\prem}\left\{L(\sym^rf,s), f\in \primeps{\kappa}{q}\right\}.%
\]
Indeed, we compute the asymptotic expectation $\Eheps[\infty]\left(D_{1}[\Phi;r]\right)$.
\begin{theorem}
\label{density2}
Let $r\geq 1$ be an odd integer, $\epsilon=\pm 1$ and $\Phi\in\Schwartz_\nu(\R)$. We assume that hypothesis $\Nice(r,f)$ holds 
for any prime number $q$ and any $f\in \prim{\kappa}{q}$ and also that $\theta$ is admissible. Let
\[%%
\nu_{1,\mathrm{max}}^\epsilon(r,\kappa,\theta)\coloneqq\inf{\left(\nu_{1,\mathrm{max}}(r,\kappa,\theta),\frac{3}{r(r+2)}\right).}%
\]
If $\nu<\nu_{1,\mathrm{max}}^\epsilon(r,\kappa,\theta)$ then
\[%%
\Eheps[\infty]\left(D_{1}[\Phi;r]\right)=E[\Phi;r].%
\]
\end{theorem}
Some comments are given in remark \ref{remark5} page \pageref{remark5}.
\begin{proof}[\proofname{} of theorem \ref{density2}] By \eqref{eq_removeps}, we have
\begin{equation}\label{eq_fmnr}%%
\Eheps[q]\left(D_{1,q}[\Phi;r]\right)%
= %
\Eh[q]\left(D_{1,q}[\Phi;r]\right)-%
\epsilon\times\epsilon(k,r)\sqrt{q}\Eh[q]\left(\lambda_.(q)D_{1,q}[\Phi;r]\right).
%\epsilon\times\epsilon(k,r)\sqrt{q}\sumh_{f\in\prim{\kappa}{q}}\lambda_f(q)D_{1,q}[\Phi;r](f).%
\end{equation}
The first term is the main term of the theorem thanks to theorem \ref{density1}. 
According to proposition \ref{explicit} %
and corollary~\ref{lem_sumlambdafq}%
, the second term (without the epsilon factors) is given by
\begin{multline}
\label{start}
\sqrt{q}\Eh[q]\left(\lambda_.(q)P_q^1[\Phi;r]\right) \\%
+\sqrt{q}\sum_{m=0}^{r-1}(-1)^m\Eh[q]\left(\lambda_.(q)P_q^2[\Phi;r,m]\right)+O\left(\frac{1}{\log{\left(q^r\right)}}\right).%
\end{multline}
Let us focus on the %
first %
term in \eqref{start} knowing that the same
discussion holds for the  %
second %
term with even 
better results on $\nu$. We have
\begin{equation}\label{eq_resun}%%
\sqrt{q}\Eh[q]\left(\lambda_.(q)P_q^1[\Phi;r]\right)=\sqrt{q}\PP{q,\mathrm{new}}1[\Phi;r]+\sqrt{q}\PP{q,\mathrm{old}}1[\Phi;r]%
\end{equation}
where
\begin{align*}
\PP{q,\mathrm{new}}1[\Phi;r]%
&=%
-\frac{2}{\log{\left(q^r\right)}}\sum_{\substack{p\in\prem \\
p\nmid q}}\Delta_q\left(p^{r}q,1\right)\frac{\log{p}}{\sqrt{p}}\widehat{\Phi}\left(\frac{\log{p}}{\log{\left(q^r\right)}}\right), \\%
\PP{q,\mathrm{old}}1[\Phi;r]%
&=%
\frac{2}{q\nu(q)\log{\left(q^r\right)}}\sum_{\ell\mid q^\infty}\frac{1}{\ell}\sum_{\substack{p\in\prem \\%
p\nmid q}}\Delta_1\left(p^{r}\ell^2,q\right)\frac{\log{p}}{\sqrt{p}}\widehat{\Phi}\left(\frac{\log{p}}{\log{\left(q^r\right)}}\right).%
\end{align*}
Lemma \ref{lem_old} implies %
\begin{equation}\label{eq_rqpuqo}%
\sqrt{q}\PP{q,\mathrm{old}}1[\Phi;r]\ll q^{(\nu r-4)/2}%
\end{equation} %
which is an admissible error term if $\nu<4/r$. %
The new part is given by
\[%
\PP{q,\mathrm{new}}1[\Phi;r]%
=%
-\frac{2%
(2\pi i^\kappa)
}{\log{\left(q^r\right)}}%
\sum_{\substack{c\geq 1 \\%
q\mid c}}\sum_{\substack{p\in\prem \\%
q\nmid p}}\frac{\log{p}}{\sqrt{p}}\frac{S\left(p^{r}q,1;c\right)}{c}%
J_{\kappa-1}\left(\frac{4\pi\sqrt{p^{r}q}}{c}\right)%
\widehat{\Phi}\left(\frac{\log{\left(p\right)}}{\log{\left(q^{r}\right)}}\right).%
\]
and can be written as
\begin{equation*}
-\frac{2(2\pi i^\kappa)}{\log{\left(q^r\right)}}\sum_{\substack{c\geq 1 \\
q\mid c}}\sum_{m\geqslant 1}\widehat{a}_m\frac{S(m,1;c)}{c}J_{\kappa-1}\left(\frac{4\pi\sqrt{m}}{c}\right)\widehat{\Phi}\left(\frac{\log{\left(m/q\right)}}{\log{(q^{r^2})}}\right)
\end{equation*}
where
\begin{equation*}
\widehat{a}_m\coloneqq\un_{[1,q^{1+\nu r^2}]}\begin{cases}
0 & \text{if $q\nmid m$ or $m\neq p^rq$ for some $p\neq q$ in $\mathcal{P}$}, \\
\frac{\log{p}}{\sqrt{p}} & \text{if $m=p^rq$ for some $p\neq q$ in $\mathcal{P}$.}
\end{cases}
\end{equation*}
Thus, if $\nu\leq 1/r^2$ then we obtain %
\[%
\PP{q,\mathrm{new}}1[\Phi;r,m]%
\ll_\epsilon%
q^\epsilon%
\sumsh_{\substack{M\leq q^{1+\nu r^2}\\ C\geq q/2}}%
\left(\frac{\sqrt{M}}{C}\right)^{\kappa-1-2\theta}%
\left(1+\sqrt{\frac{M}{q}}\right)%
\]
as in the proof of corollary~\ref{usefulsieve}. Summing over $C$ \emph{via} \eqref{dyadic2} gives
\[%
\PP{q,\mathrm{new}}1[\Phi;r%
,m%
]\ll_\epsilon %
q^\epsilon\sumsh_{M\leq q^{1+r^2\nu}}%
\left(\frac{M^{\frac{\kappa-1}{2}-\theta}}{q^{\kappa-1-2\theta}}+%
\frac{M^{\frac{\kappa}{2}-\theta}}{q^{\kappa-\frac{1}{2}-2\theta}}\right).
\]
Summing over $M$ \emph{via} \eqref{dyadic1} leads to
\begin{equation}\label{eq_rqpuqn}
\PP{q,\mathrm{new}}1[\Phi;r%
,m%
]\ll_\epsilon %
q^{\left(\frac{\kappa-1}{2}-\theta\right)r^2\nu-(\frac{\kappa-1}{2}-\theta)+\epsilon}+%
q^{\left(\frac{\kappa}{2}-\theta\right)r^2\nu-\left(\frac{\kappa-1}{2}-\theta\right)+\epsilon}
\end{equation}
which is an admissible error term if $\nu<\frac{1}{r^2}\left(1-\frac{1}{\kappa-2\theta}\right)$.
\end{proof}
%................................
\section{Quadratic statistics for low-lying zeros}%
\label{two}%
%................................
\subsection{Asymptotic expectation of the two-level density and asymptotic variance}%
\label{sec_twoandvar}%
%................................
\begin{definition}\label{def_tld}
Let $f\in\prim{\kappa}{q}$ and $\Phi_1$, $\Phi_2$ in $\Schwartz_\nu(\R)$. The \emph{two-level density} %
(relatively to $\Phi_1$ and $\Phi_2$) of $\sym^rf$ is
\[%
D_{2,q}[\Phi_1,\Phi_2;r](f)%
\coloneqq%
\sum_{\substack{(j_1,j_2)\in\mathcal{E}(f,r)^2\\ j_1\neq\pm j_2}}%
\Phi_1\left(\widehat{\rho}_{f,r}^{(j_1)}\right)
\Phi_2\left(\widehat{\rho}_{f,r}^{(j_2)}\right).
\]
\end{definition}
\begin{remark}
In this definition, it is important to note that the condition $j_1\neq j_2$ \emph{does not imply} that %
 $\widehat{\rho}_{f,r}^{(j_1)}\neq \widehat{\rho}_{f,r}^{(j_2)}$. It only implies this if the zeros are simple. %
Recall however that some $L$-functions of elliptic curves (hence of modular forms) have multiple zeros at the critical point \cite{MR848380, %
MR777279}.
\end{remark}
The following lemma is an immediate consequence of definition~\ref{def_tld}.
\begin{lemma}\label{lem_tld}
Let $f\in\prim{\kappa}{q}$ and $\Phi_1$, $\Phi_2$ in $\Schwartz_\nu(\R)$.
Then, %
\begin{multline*}
D_{2,q}[\Phi_1,\Phi_2;r](f)=%
D_{1,q}[\Phi_1;r](f)D_{1,q}[\Phi_2;r](f)%
-2D_{1,q}[\Phi_1\Phi_2;r](f)\\
+%
%\delta\left(\epsilon\left(\sym^rf\right)=-1\right)%
\un_{\primimpair{\kappa}{q}}(f)%
\times\Phi_1(0)\Phi_2(0).
\end{multline*}
\end{lemma}
We first evaluate the product of one-level statistics on average.
\begin{lemma}\label{lem_conesson}%
Let $r\geq 1$. %
Let $\Phi_1$ and $\Phi_2$ in $\Schwartz_\nu(\R)$. %
We assume that hypothesis %
$\Nice(r,f)$ holds for any prime number $q$ %
and any $f\in \prim{\kappa}{q}$ and also that $\theta$ is admissible. %
If $\;\nu<1/r^2$ then %
\[%
\Eh[\infty]\left(%
D_{1}[\Phi_1;r]D_{1}[\Phi_2;r]%
\right)%
=%
E[\Phi_1;r]E[\Phi_2;r]+%
2\int_\R\abs{u}\widehat{\Phi_1}(u)\widehat{\Phi_2}(u)\dd u.
\]
\end{lemma}
\begin{remark}
Since theorem~\ref{density1} %
implies that %
\begin{multline*}
\Eh[\infty]\left(%
D_{1}[\Phi_1;r]D_{1}[\Phi_2;r]%
\right)%
-%
E[\Phi_1;r]E[\Phi_2;r]
=\\
\Eh[\infty]\left(%
D_{1}[\Phi_1;r]D_{1}[\Phi_2;r]%
\right)
-
\Eh[\infty]\left(D_{1}[\Phi_1;r]\right)%
\Eh[\infty]\left(D_{1}[\Phi_2;r]\right),
\end{multline*}
lemma \ref{lem_conesson} reveals %
that the term
\[%
\Ch[\infty]\left(D_{1}[\Phi_1;r],D_{1}[\Phi_2;r]\right)\coloneqq %
2\int_\R\abs{u}\widehat{\Phi_1}(u)\widehat{\Phi_2}(u)\dd u%
\]
measures the dependence between %the random variables %
$D_{1}[\Phi_1;r]$ and $D_{1}[\Phi_2;r]$. This term is the
\emph{asymptotic covariance} of $D_{1}[\Phi_1;r]$ and $D_{1}[\Phi_2;r]$. In
particular, taking $\Phi_1=\Phi_2$, we obtain the asymptotic variance.
\end{remark}
\begin{theorem}
\label{th_variance}
Let $\Phi\in\Schwartz_{\nu}(\R)$. If $\nu<1/r^2$ then the asymptotic variance of %
the random variable %
$D_{1,q}[\Phi;r]$ is %
\[%
\Vh[\infty]\left(D_1[\Phi;r]\right)=2\int_\R\abs{u}\widehat{\Phi}^2(u)\dd u.
\]
\end{theorem}
\begin{proof}[\proofname{} of lemma~\ref{lem_conesson}]
From proposition~\ref{explicit}, we obtain %
\begin{multline}
\label{eq_prodD}
\Eh[q]\left(D_{1,q}[\Phi_1;r]D_{1,q}[\Phi_2;r]\right)%
=%
E[\Phi_1;r]E[\Phi_2;r]+\Ch[q]%
\\
+\sum_{\substack{(i,j)\in\{1,2\}^2\\ i\neq j}}%
\sum_{m=0}^{r-1}(-1)^m\Eh[q]\left(P_q^1[\Phi_i;r]P_q^2[\Phi_j;r,m]\right)%
\\
+\sum_{m_1=0}^{r-1}\sum_{m_2=0}^{r-1}(-1)^{m_1+m_2}\Eh[q]\left(P_q^2[\Phi_1;r,m_1]P_q^2[\Phi_2;r,m_2]\right)%
+O\left(\frac{1}{\log{\left(q^r\right)}}\right)
\end{multline}
with
\[%
\Ch[q]\coloneqq \Eh[q]\left(P_q^1[\Phi_1;r]P_q^1[\Phi_2;r]\right).
\]
The error term is evaluated by use of theorem~\ref{density1} and
equations~\eqref{eq_mesasy}, \eqref{eq_psum} and \eqref{eq_psdm}.
We first compute $\Ch[q]$. Using proposition~\ref{iwlusatr}, we compute
$\Ch[q]=E^n-4E^o$ with
\begin{equation*}
E^n\coloneqq%
\frac{4}{\log^2{\left(q^r\right)}}%
\sum_{\substack{p_1\in\prem\\ p_1\nmid q}}%
\sum_{\substack{p_2\in\prem\\ p_2\nmid q}}%
\frac{\log p_1}{\sqrt{p_1}}\frac{\log p_2}{\sqrt{p_2}}%
\widehat{\Phi_1}\left(\frac{\log p_1}{\log{\left(q^r\right)}}\right)%
\widehat{\Phi_2}\left(\frac{\log p_2}{\log{\left(q^r\right)}}\right)%
\Delta_q(p_1^r,p_2^r)
\end{equation*}
and %
\begin{multline*}
E^o\coloneqq%
\frac{1}{q\log^2{\left(q^r\right)}}%
\\%
\times%
\sum_{\substack{p_1\in\prem\\ p_1\nmid q}}%
\sum_{\substack{p_2\in\prem\\ p_2\nmid q}}%
\frac{\log p_1}{\sqrt{p_1}}\frac{\log p_2}{\sqrt{p_2}}%
\widehat{\Phi_1}\left(\frac{\log p_1}{\log{\left(q^r\right)}}\right)%
\widehat{\Phi_2}\left(\frac{\log p_2}{\log{\left(q^r\right)}}\right)%
\sum_{\ell\mid q^\infty}\frac{\Delta_1(\ell^2p_1^r,p_2^r)}{\ell}.
\end{multline*}
By definition of the $\Delta$-symbol, we write
$E^n=E^n_{\mathrm{p}}+\frac{8\pi i^{\kappa}}{\log^2{\left(q^r\right)}}E^n_{\mathrm{e}}$ with
\[%
E^n_{\mathrm{p}}\coloneqq%
\frac{4}{\log^2{(q^r)}}%
\sum_{\substack{p\in\prem\\ p\nmid q}}%
\frac{\log^2 p}{p}\left(\widehat{\Phi_1}\widehat{\Phi_2}\right)%
\left(\frac{\log p}{\log{(q^r)}}\right)%
\]
and
\begin{multline*}
E^n_{\mathrm{e}}\coloneqq%
\sum_{\substack{c\geq 1\\ q\mid c}}
\sum_{\substack{p_1\in\prem\\ p_1\nmid q}}%
\sum_{\substack{p_2\in\prem\\ p_2\nmid q}}%
\frac{\log p_1}{\sqrt{p_1}}\frac{\log p_2}{\sqrt{p_2}}%
\widehat{\Phi_1}\left(\frac{\log p_1}{\log{\left(q^r\right)}}\right)%
\widehat{\Phi_2}\left(\frac{\log p_2}{\log{\left(q^r\right)}}\right)%
\\
\times%
\frac{S(p_1^r,p_2^r;c)}{c}%
J_{\kappa-1}\left(\frac{4\pi\sqrt{p_1^rp_2^r}}{c}\right).
\end{multline*}
We remove the condition $p\nmid q$ from $E^n_{\mathrm{p}}$ at an admissible cost
and obtain, after integration by parts,
\begin{equation}\label{eq_epn}
E^n_{\mathrm{p}}=2\int_\R\abs{u}\widehat{\Phi_1}(u)\widehat{\Phi_2}(u)\dd u%
+O\left(\frac{1}{\log^2{(q^r)}}\right).
\end{equation}%
%.............
Using corollary~\ref{usefulsieve}, we get
\begin{equation}\label{eq_een}
E^n_{\mathrm{e}}\ll\frac{1}{\log^2{\left(q^r\right)}}
\end{equation}
as soon as $\nu\leq 1/r^2$.
%.............
Finally, using lemma~\ref{lem_old}, we see that $E^o$ is an admissible error term for $\nu<1/r$ so that
equations \eqref{eq_epn} and \eqref{eq_een} lead to
\begin{equation}\label{eq_en}
\Ch[q]=2\int_\R\abs{u}\widehat{\Phi_1}(u)\widehat{\Phi_2}(u)\dd u%
+O\left(\frac{1}{\log^2{\left(q^r\right)}}\right).
\end{equation}
Let $\{i,j\}=\{1,2\}$.
We prove next that each $\Eh[q]\left(P^1_q[\Phi_i;r]P^2_q[\Phi_j;r,m]\right)$ is
an error term when $\nu<1/r^2$. %
Using proposition~\ref{iwlusatr} and %
lemma~\ref{lem_old} we have %
\begin{multline*}
\Eh[q]\left(P^1_q[\Phi_i;r]P^2_q[\Phi_j;r,m]\right)%
=%
\frac{8\pi i^\kappa}{\log^2{(q^r)}}%
\sum_{\substack{c\geq 1\\ q\mid c}}%
\sum_{\substack{p_1\in\prem\\ p_1\nmid q}}%
\sum_{\substack{p_2\in\prem\\ p_2\nmid q}}%
\frac{\log p_1}{\sqrt{p_1}}\frac{\log p_2}{p_2}%
\widehat{\Phi_i}\left(\frac{\log p_1}{\log{\left(q^r\right)}}\right)%
\\
\times%
\widehat{\Phi_j}\left(\frac{\log p_2}{\log{\left(q^{r/2}\right)}}\right)%
\frac{S(p_1^r,p_2^{2r-2m};c)}{c}%
J_{\kappa-1}\left(\frac{4\pi\sqrt{p_1^rp_2^{2r-2m}}}{c}\right)
+O\left(\frac{1}{\log{\left(q^r\right)}}\right)^2.
\end{multline*}
%............
We use corollary~\ref{usefulsieve} to conclude that%
\begin{equation}\label{eq_utile}
\Eh[q]\left(P^1_q[\Phi_i;r]P^2_q[\Phi_j;r,m]\right)
\ll\frac{1}{\log q}
\end{equation}
when $\nu<1/r^2$. %
%............
Finally, %
$\Eh[q]\left(P^2_q[\Phi_1;r,m_1]P^2_q[\Phi_2;r,m_2]\right)$ is shown to %
be an error term in the same way.
\end{proof}
Using lemmas~\ref{lem_tld} and \ref{lem_conesson}, theorem~\ref{density1},
hypothesis~$\Nice(r,f)$ and remark~\ref{rem_moyun}, we prove the following
theorem.
\begin{theorem}\label{thm_twodens}%
%.............

%%
Let $r\geq 1$. %
Let $\Phi_1$ and $\Phi_2$ in $\Schwartz_\nu(\R)$. %
We assume that hypothesis %
$\Nice(r,f)$ holds for any prime number $q$ %
and any $f\in \prim{\kappa}{q}$ and also that $\theta$ is admissible. %
If $\;\nu<\nu_{2,\mathrm{max}}(r,\kappa,\theta)$ then %

%%
%.............
\begin{multline*}
\Eh[\infty]\left(D_2[\Phi_1,\Phi_2;r]\right)%
=%
\left[%
\widehat{\Phi_1}(0)%
+%
\frac{(-1)^{r+1}}{2}\Phi_1(0)%
\right]%
\left[%
\widehat{\Phi_2}(0)%
+%
\frac{(-1)^{r+1}}{2}\Phi_2(0)%
\right]%
\\
+2\int_\R\abs{u}\widehat{\Phi_1}(u)\widehat{\Phi_2}(u)\dd u
-2\widehat{\Phi_1\Phi_2}(0)%
+\left((-1)^{r}+%
\frac{%
%\delta(2\nmid r)%
\un_{2\N+1}(r)
}{2}\right)%
\Phi_1(0)\Phi_2(0).%
\end{multline*}
\end{theorem} 
Some comments are given in remark \ref{remark6} page \pageref{remark6}.
%................................
\subsection{Signed asymptotic expectation of the two-level density and signed asymptotic variance}
In this part, $r$ is \emph{odd}. 

\begin{lemma}\label{lem_signedconesson}
Let $\Phi_1$ and $\Phi_2$ in $\Schwartz_\nu(\R)$.
If $\nu<1/(2r^2)$ then %
\[%
\Eheps[\infty]\left(%
D_{1}[\Phi_1;r]D_{1}[\Phi_2;r]%
\right)%
=%
E[\Phi_1;r]E[\Phi_2;r]+%
2\int_\R\abs{u}\widehat{\Phi_1}(u)\widehat{\Phi_2}(u)\dd u.
\]
\end{lemma}

\begin{remark}
By  theorem~\ref{density2} and lemma~\ref{lem_signedconesson} we have
\begin{multline*}
\Eheps[\infty]\left(%
D_{1}[\Phi_1;r]D_{1}[\Phi_2;r]%
\right)%
-%
E[\Phi_1;r]E[\Phi_2;r]
=\\
\Eheps[\infty]\left(%
D_{1}[\Phi_1;r]D_{1}[\Phi_2;r]%
\right)
-
\Eheps[\infty]\left(D_{1}[\Phi_1;r]\right)%
\Eheps[\infty]\left(D_{1}[\Phi_2;r]\right).
\end{multline*}
Thus,
\[%
\Cheps[\infty]\left(D_{1}[\Phi_1;r],D_{1}[\Phi_2;r]\right)\coloneqq %
2\int_\R\abs{u}\widehat{\Phi_1}(u)\widehat{\Phi_2}(u)\dd u
\]
is the \emph{signed asymptotic covariance} of $D_{1}[\Phi_1;r]$ and $D_{1}[\Phi_2;r]$. In
particular, taking $\Phi_1=\Phi_2$, we obtain the signed asymptotic variance.
\end{remark}
\begin{theorem}
Let $\Phi\in\Schwartz_{\nu}(\R)$. If $\nu<1/(2r^2)$ then the signed asymptotic variance of
$D_1[\Phi;r]$ is
\[%
\Vheps[\infty]\left(D_1[\Phi;r]\right)=2\int_\R\abs{u}\widehat{\Phi}^2(u)\dd u.
\]
\end{theorem}
\begin{proof}[\proofname{} of lemma~\ref{lem_signedconesson}]%
From proposition~\ref{explicit} and \eqref{eq_mbun}, we obtain %
\begin{multline}
\label{eq_resu}
\Eheps[q]\left(D_{1,q}[\Phi_1;r]D_{1,q}[\Phi_2;r]\right)%
=%
E[\Phi_1;r]E[\Phi_2;r]%
+\Cheps[q]%
\\%
+%
\sum_{\substack{(i,j)\in\{1,2\}^2\\ i\neq j}}%
\sum_{m=0}^{r-1}(-1)^m\Eheps[q]\left(P_q^1[\Phi_i;r]P_q^2[\Phi_j;r,m]\right)%
\\
+\sum_{m_1=0}^{r-1}\sum_{m_2=0}^{r-1}(-1)^{m_1+m_2}\Eheps[q]\left(P_q^2[\Phi_1;r,m_1]P_q^2[\Phi_2;r,m_2]\right)%
+O\left(\frac{1}{\log{\left(q^r\right)}}\right)
\end{multline}
with
\[%
\Cheps[q]\coloneqq \Eheps[q]\left(P_q^1[\Phi_1;r]P_q^1[\Phi_2;r]\right).
\]%
Assume that $\nu<1/r^2$. Then equations~\eqref{eq_removeps}, \eqref{eq_en} %
and proposition~\ref{iwlusatr} lead to %
\begin{equation}\label{eq_valCheps}%
\Cheps[q]=2\int_{\R}\abs{u}\widehat{\Phi_1}(u)\widehat{\Phi_2}(u)\dd u%
-\epsilon\times\epsilon(\kappa,r)(G^n-4G^o)
\end{equation}
with
\[%%
G^n\coloneqq \frac{4\sqrt{q}}{\log^2{\left(q^r\right)}}%
\sum_{\substack{p_1\in\prem\\ p_1\nmid q}}%
\sum_{\substack{p_2\in\prem\\ p_2\nmid q}}%
\frac{\log p_1}{\sqrt{p_1}}%
\frac{\log p_2}{\sqrt{p_2}}%
\widehat{\Phi_1}\left(\frac{\log p_1}{\log{\left(q^r\right)}}\right)%
\widehat{\Phi_2}\left(\frac{\log p_2}{\log{\left(q^r\right)}}\right)%
\Delta_q\left(p_1^rq,p_2^r\right)
\]
and
\begin{multline}%
G^o\coloneqq \frac{1}{\sqrt{q}\log^2{\left(q^r\right)}}%
\\%
\times%
\sum_{\substack{p_1\in\prem\\ p_1\nmid q}}%
\sum_{\substack{p_2\in\prem\\ p_2\nmid q}}%
\frac{\log p_1}{\sqrt{p_1}}%
\frac{\log p_2}{\sqrt{p_2}}%
\widehat{\Phi_1}\left(\frac{\log p_1}{\log{\left(q^r\right)}}\right)%
\widehat{\Phi_2}\left(\frac{\log p_2}{\log{\left(q^r\right)}}\right)%
\sum_{\ell\mid q^\infty}%
\frac{\Delta_q\left(\ell^2p_1^r,p_2^rq\right)}{\ell}.
\end{multline}
Lemma~\ref{deltaestimate} implies that if $\nu<1/(2r^2)$ then %
\begin{equation}\label{eq_majoGn}
G^n\ll
\frac{q^{\nu r[r(\kappa-1)+1]/2}}{q^{(\kappa-1)/2}}
%\frac{q^{-1/4+2\nu r(r/4+1/2)}}{\log^2{\left(q^r\right)}}
\end{equation}
hence $G^n$ is an error term as soon as $\nu\leq 1/(2r^2)$.
Lemma~\ref{lem_old} implies %
\begin{equation}\label{eq_majoGo}
G^o\ll q^{-3/2+\nu r+\epsilon}
\end{equation}
which is an error term. %
Reporting equations~\eqref{eq_majoGn} and \eqref{eq_majoGo} in %
\eqref{eq_valCheps} %
we obtain %
\begin{equation}\label{eq_unos}
\Cheps[\infty]=%
2\int_{\mathbb{R}}\abs{u}\widehat{\Phi_1}(u)\widehat{\Phi_2}(u)\dd u%
\end{equation}
for $\nu\leq 1/(2r(r+2))$. Next, we prove that each %
$\Eheps[q]\left(P_q^1[\Phi_i;r]P_q^2[\Phi_j;r,m]\right)$ %
is an error term as soon as %
$\nu\leq 1/(2r^2)$. From equations~\eqref{eq_removeps} and \eqref{eq_utile}, %
we obtain %
\begin{multline}\label{eq_signedsecond}
\Eheps[q]\left(P_q^1[\Phi_i;r]P_q^2[\Phi_j;r,m]\right)%
=\\
-\epsilon\times\epsilon(\kappa,r)\sqrt{q}\sumh_{f\in\prim{\kappa}{q}}%
\lambda_f(q)P_q^1[\Phi_i;r]P_q^2[\Phi_j;r,m]%
+O\left(\frac{1}{\log q)}\right).%
\end{multline}
We use proposition~\ref{iwlusatr} and lemmas~\ref{lem_old} %
and \ref{deltaestimate} to have %
\begin{multline}\label{eq_signedinter}
\sqrt{q}\sumh_{f\in\prim{\kappa}{q}}%
\lambda_f(q)P_q^1[\Phi_i;r]P_q^2[\Phi_j;r,m]%
\ll\\%
\frac{q^{\nu r(2r-m+2)/4-1/4}}{\log^2q}+%
\frac{q^{(\nu r-1)/2+\epsilon}}{\log q}.%
\end{multline}
It follows from ~\eqref{eq_signedinter} and \eqref{eq_signedsecond} that %
\begin{equation}\label{eq_dos}
\Eheps[\infty]\left(P_q^1[\Phi_i;r]P_q^2[\Phi_j;r,m]\right)=0
\end{equation}
for $\nu\leq 1/(2r(r+1))$. In the same way, we have, for $\nu$ in the previous range,
\begin{equation}\label{eq_tres}
\Eheps[\infty]\left(P_q^2[\Phi_1;r,m_1]P_q^2[\Phi_2;r,m_2]\right)=0.
\end{equation}
Reporting \eqref{eq_unos}, \eqref{eq_dos} and \eqref{eq_tres} in \eqref{eq_resu}, we have the announced result. 
\end{proof}

Using lemmas~\ref{lem_tld}, \ref{lem_signedconesson}, theorem~\ref{density2},
hypothesis~$\Nice(r,f)$ and \eqref{eq_mbun}, we prove the following
theorem.

\begin{theorem}\label{thm_signedtwodens}
Let $f\in\prim{\kappa}{q}$ and $\Phi_1$, $\Phi_2$ in $\Schwartz_\nu(\R)$.
If $\nu<1/(2r(r+1))$ then
\begin{multline*}
\Eheps[\infty]\left(D_2[\Phi_1,\Phi_2;r]\right)%
=%
\left[%
\widehat{\Phi_1}(0)%
+%
\frac{1}{2}\Phi_1(0)%
\right]%
\left[%
\widehat{\Phi_2}(0)%
+%
\frac{1}{2}\Phi_2(0)%
\right]%
\\
+2\int_\R\abs{u}\widehat{\Phi_1}(u)\widehat{\Phi_2}(u)\dd u
-2\widehat{\Phi_1\Phi_2}(0)%
-\Phi_1(0)\Phi_2(0)%
\\
+%
%\delta(\epsilon=-1)%
\un_{\{-1\}}(\epsilon)
\Phi_1(0)\Phi_2(0).
\end{multline*}
\end{theorem} 

\begin{remark}
Remark~\ref{rem_symtyp} together with theorem~\ref{thm_signedtwodens} and a result of Katz \& Sarnak %
(see \cite[Theorem A.D.2.2]{KaSa} or \cite[Theorem 3.2]{Mil}) imply that the symmetry type of $\mathcal{F}_r^{\epsilon}$
is as in table~\ref{tab_symty}. Some additional comments are given in remark \ref{remark2} page \pageref{remark2}.
\begin{table}[ht]
%\begin{center}
\setlength{\extrarowheight}{4pt}
\begin{tabular}{|c|c|c|}
\hline
\backslashbox{$\epsilon$}{$r$} & even & odd\\ 
\hline
$-1$ &  & $SO(\mathrm{odd})$ \\
\hline
$1$ & $Sp$  & $SO(\mathrm{even})$\\
\hline
\end{tabular}
%\end{center}
\caption{Symmetry type of $\mathcal{F}_r^{\epsilon}$}\label{tab_symty}
\end{table}
\end{remark}
%................................
%\subsection{Decomposition of the variance}
\section{First asymptotic moments of the one-level density}

\label{momentt}

In this section, we compute the asymptotic $m$-th moment of the one level density namely
\begin{equation*}
\Mh[\infty,m]\left(D_{1,q}[\Phi;r]\right)\coloneqq\lim_{\substack{q\in\mathcal{P} \\ q\to+\infty}}\Mh[q,m]\left(D_{1,q}[\Phi;r]\right)
\end{equation*}
where
\begin{equation*}
\Mh[q,m]\left(D_{1,q}[\Phi;r]\right)=%
\Eh[q]%
\left(%
\left(D_{1,q}[\Phi;r]-\Eh[q](D_{1,q}[\Phi;r])\right)^m
\right)%
\end{equation*}
for $m$ small enough (regarding to the size of the support of $\Phi$). 
The end of this section is devoted to the proof of theorem \ref{thm_F}. %
Note that we can assume that $m\geq 3$ since the work has already been done
for $m=1$ and $m=2$. %
Thanks to equation~\eqref{eq_mhe} and proposition~\ref{explicit}, we have
\begin{align}
\Mh[q,m]\left(D_{1,q}[\Phi;r]\right)  & = %
 \sum_{\ell=0}^m\binom{m}{\ell}\Eh[q]\left(P_q^1[\Phi;r]^{m-\ell}\left(P_q^2[\Phi;r]+%
O\left(\frac{1}{\log q}\right)\right)^\ell\right) \\
\label{eq_evmom} \\
& =  \sum_{\substack{0\leq\ell\leq m \\
0\leq\alpha\leq\ell}}\binom{m}{\ell}\binom{\ell}{\alpha}R(q)^{\ell-\alpha}%
\Eh[q]\left(P_q^1[\Phi;r]^{m-\ell}P_q^2[\Phi;r]^{\alpha}\right)
\end{align}
where % 
\begin{align}
P_q^2[\Phi;r](f) & \coloneqq %
 -\frac{2}{\log(q^r)}%
\sum_{j=0}^{r-1}(-1)^j%
\sum_{\substack{p\in\prem\\ p\nmid q}}%
\lambda_f\left(p^{2(r-j)}\right)\frac{\log p}{p}\widehat{\Phi}\left(\frac{2\log p}{\log(q^r)}\right) \\
\label{eq_tard} & = %
 -\frac{2}{\log(q^r)}\sum_{j=1}^{r}(-1)^{r-j}\sum_{\substack{p\in\prem\\ p\nmid %
  q}}\lambda_f\left(p^{2j}\right)\frac{\log p}{p}\widehat{\Phi}\left(\frac{2\log %
   p}{\log(q^r)}\right)
\end{align}
and $R$ is a positive function satisfying %
\[%
R(q)\ll\frac{1}{\log q}.
\]
Thus, an asymptotic formula for $\Mh[q,m]\left(D_{1,q}[\Phi;r]\right)$ directly follows from the next proposition.
\begin{proposition}
\label{moments}
Let $r\geq 1$ be any integer. We assume that hypothesis $\Nice(r,f)$ %
holds for any prime number $q$ and any primitive holomorphic cusp form of level $q$ and even weight $\kappa$. Let $\alpha\geq 0$ and $\ell\geq 0$ be any integers.
\begin{itemize}
\item
If $\;\alpha\geq 1$ and $\;\alpha\nu<4/r^2$ then
\begin{equation*}
\Eh[q]\left(P_q^2[\Phi;r]^{\alpha}\right)=O\left(\frac{1}{\log{q}}\right).
\end{equation*}
\item
If $\;1\leq\alpha\leq\ell\leq m-1$ and $(\alpha+m-\ell)\nu<4/(r(r+2))$ then
\begin{equation*}
\Eh[q]\left(P_q^1[\Phi;r]^{m-\ell}P_q^2[\Phi;r]^{\alpha}\right)=O\left(\frac{1}{\log{q}}\right).
\end{equation*}
\item
If $\;\alpha\geq 1$ and $\;\alpha\nu<4/(r(r+2))$ then
\begin{equation*}
\Eh[q]\left(P_q^1[\Phi;r]^{\alpha}\right)=\begin{cases}
O\left(\frac{1}{\log^2{(q)}}\right) & \text{if $\alpha$ is odd,} \\
2\int_\R\abs{u}\widehat{\Phi}^2(u)\dd u\times\frac{\alpha!}{2^{\alpha/2}\left(\frac{\alpha}{2}\right)!}+O\left(\frac{1}{\log^2{(q)}}\right) & \text{otherwise}.
\end{cases}
\end{equation*}
\end{itemize}
\end{proposition}
\subsection{One some useful combinatorial identity} In order to use the multiplicative properties of Hecke eigenvalues in the proof of proposition \ref{moments}, we want to reorder some sums over many primes to sums
over distinct primes. We follow the work of Hughes \& Rudnick \cite[\S 7]{HuRu} (see also %
\cite{MR2166468} and the work of Soshnikov \cite{Sos00}) to achieve this. %
Let $P(\alpha,s)$ be the set of surjective functions %
\[%
\sigma \colon \{1,\dotsc,\alpha\}\twoheadrightarrow \{1,\dotsc,s\}%
\]
such that for any $j\in\{1,\dotsc,\alpha\}$, either $\sigma(j)=1$ or there exists $k<j$ such that $\sigma(j)=\sigma(k)+1$.
This can be viewed as the number of partitions of a set of $\alpha$ elements into $s$ nonempty subsets.
By definition, the cardinality of $P(\alpha,s)$ is the Stirling number of second kind \cite[\S 1.4]{Sta97}.
For any $j\in\{1,\dotsc,s\}$, let %
\[%
\varpi_j^{(\sigma)}\coloneqq\#\sigma^{-1}(\{j\}).%
\]
Note that %
\begin{equation}\label{eq_trivcond}%
\varpi_j^{(\sigma)}\geq 1 \quad\text{ for any $1\leq j\leq s$}\qquad\text{ and }\qquad
\sum_{j=1}^s\varpi_j^{(\sigma)}=\alpha.
\end{equation}
The following lemma is lemma 7.3 of \cite[\S 7]{HuRu}.
\begin{lemma}
\label{primesdistincts}
If $g$ is any function of $m$ variables then
\begin{equation*}
\sum_{j_1,\dotsc,j_m}g\left(x_{j_1},\dotsc,x_{j_m}\right)=\sum_{s=1}^m\sum_{\sigma\in P(m,s)}\sum_{\substack{i_1,\dotsc,i_s \\
\text{distinct}}}g\left(x_{i_{\sigma(1)}},\dotsc,x_{i_{\sigma(m)}}\right).
\end{equation*}
\end{lemma}
\subsection{Proof of the first bullet of proposition \ref{moments}} By the definition \eqref{eq_tard}, we have
\begin{multline}\label{eq_mais}
\Eh[q]\left(P_q^2[\Phi;r]^\alpha\right)=\frac{(-2)^\alpha}{\log^\alpha{(q^r)}}\sum_{1\leq j_1,\dotsc,j_\alpha\leq r}(-1)^{\alpha r-(j_1+\dotsc+j_\alpha)} \\
\times\sum_{\substack{p_1,\dotsc,p_\alpha\in\prem\\ q\nmid p_1\dotsc p_\alpha}}\left(\prod_{i=1}^\alpha\frac{\log p_i}{p_i}\widehat{\Phi}\left(\frac{2\log p_i}{\log(q^r)}\right)\right)\Eh[q]\left(\prod_{i=1}^\alpha\lambda_f\left(p_i^{2j_i}\right)\right).
\end{multline}
Writing $\{\widehat{p}_i\}_{i\geq 1}$ for the increasing sequence of prime numbers except $q$, we have
\begin{multline}\label{eq_renum}
\sum_{\substack{p_1,\dotsc,p_\alpha\in\prem\\ q\nmid p_1\dotsc p_\alpha}}\left(\prod_{i=1}^\alpha\frac{\log p_i}{p_i}\widehat{\Phi}\left(\frac{2\log p_i}{\log(q^r)}\right)\right)\Eh[q]\left(\prod_{i=1}^\alpha\lambda_f\left(p_i^{2j_i}\right)\right) \\
=\sum_{i_1,\dotsc,i_\alpha}\left(\prod_{\ell=1}^\alpha\frac{\log\widehat{p}_{i_\ell}}{\widehat{p}_{i_\ell}}\widehat{\Phi}\left(\frac{2\log\widehat{p}_{i_\ell}}{\log(q^r)}\right)\right)\Eh[q]\left(\prod_{\ell=1}^\alpha\lambda_f\left(\widehat{p}_{i_\ell}^{2j_\ell}\right)\right).
\end{multline}
Using lemma \ref{primesdistincts}, we rewrite the right sum in \eqref{eq_renum} as
\begin{multline}\label{eq_sosru}
\sum_{s=1}^\alpha\sum_{\sigma\in P(\alpha,s)}\sum_{\substack{k_1,\dotsc,k_s\\\text{distinct}}}\left(\prod_{i=1}^\alpha\frac{\log\widehat{p}_{k_{\sigma(i)}}}{\widehat{p}_{k_{\sigma(i)}}}\widehat{\Phi}\left(\frac{
2\log\widehat{p}_{k_{\sigma(i)}}}{\log(q^r)}\right)\right)\Eh[q]\left(\prod_{i=1}^\alpha\lambda_f\left(\widehat{p}_{k_{\sigma(i)}}^{2j_i}\right)\right) \\
=\sum_{s=1}^\alpha\sum_{\sigma\in P(\alpha,s)}\sum_{\substack{k_1,\dotsc,k_s\\\text{distinct}}}\left(\prod_{u=1}^s\left(\frac{\log\widehat{p}_{k_{u}}}{\widehat{p}_{k_{u}}}\widehat{\Phi}\left(\frac{
2\log\widehat{p}_{k_{u}}}{\log(q^r)}\right)\right)^{\varpi^{(\sigma)}_u}\right)\Eh[q]\left(\prod_{\substack{1\leq u\leq s \\
1\leq j\leq r}}\lambda_f\left(\widehat{p}_{k_{u}}^{2j}\right)^{\varpi_{u,j}^{(\sigma)}}\right)
\end{multline}
where
\begin{equation*}
\varpi_{u,j}^{(\sigma)}\coloneqq\#\{1\leq i\leq\alpha, \sigma(i)=u, j_i=j\}
\end{equation*}
for any $1\leq u\leq s$ and any $1\leq j\leq r$. Now, we show that
\begin{multline}\label{s<alpha}
\sum_{s=1}^{\alpha-1}\sum_{\sigma\in P(\alpha,s)}\sum_{\substack{k_1,\dotsc,k_s\\\text{distinct}}}\left(\prod_{u=1}^s\left(\frac{\log\widehat{p}_{k_{u}}}{\widehat{p}_{k_{u}}}\widehat{\Phi}\left(\frac{
2\log\widehat{p}_{k_{u}}}{\log(q^r)}\right)\right)^{\varpi^{(\sigma)}_u}\right)\Eh[q]\left(\prod_{\substack{1\leq u\leq s \\
1\leq j\leq r}}\lambda_f\left(\widehat{p}_{k_{u}}^{2j}\right)^{\varpi_{u,j}^{(\sigma)}}\right) \\
\ll\log^{\alpha-1}{(q)}.
\end{multline}
For $s<\alpha$ and $\sigma\in P(\alpha,s)$, we use \eqref{individualh} together with \eqref{eq_moyun} %
to obtain that the left-hand side of the previous equation is bounded by
\begin{equation}
\sum_{s=1}^{\alpha-1}\sum_{\sigma\in P(\alpha,s)}%
\sum_{\substack{k_1,\dotsc,k_s\\\text{distinct}}}%
\prod_{u=1}^s%
\left(%
\frac{%
\log\widehat{p}_{k_{u}}%
}{%
\widehat{p}_{k_{u}}%
}%
\abs{%
\widehat{\Phi}\left(%
\frac{%
2\log\widehat{p}_{k_{u}}%
}{\log(q^r)}%
\right)%
}%
\right)^{\varpi_u^{(\sigma)}}.%
\end{equation}
Since $s<\alpha$, equation~\eqref{eq_trivcond} implies that $\varpi_u^{(\sigma)}>1$ for some $1\leq u\leq s$. %
These values lead to convergent, hence bounded, sums. Let
\[%
d^{(\sigma)}\coloneqq\#\left\{1\leq u\leq s \colon \varpi_u^{(\sigma)}=1\right\}\in\{0,\dotsc,\alpha-1\},%
\]
then
\begin{multline}
\sum_{s=1}^{\alpha-1}\sum_{\sigma\in P(\alpha,s)}%
\sum_{\substack{k_1,\dotsc,k_s\\\text{distinct}}}%
\prod_{u=1}^s%
\left(%
\frac{%
\log\widehat{p}_{k_{u}}%
}{%
\widehat{p}_{k_{u}}%
}%
\abs{%
\widehat{\Phi}\left(%
\frac{%
2\log\widehat{p}_{k_{u}}%
}{\log(q^r)}%
\right)%
}%
\right)^{\varpi_u^{(\sigma)}} \\
\ll%
\sum_{s=1}^{\alpha-1}\sum_{\sigma\in P(\alpha,s)}%
\sum_{\substack{k_1,\dotsc,k_d\\\text{distinct}}}%
\prod_{u=1}^{d^{(\sigma)}}%
\left(%
\frac{%
\log\widehat{p}_{k_{u}}%
}{%
\widehat{p}_{k_{u}}%
}%
\abs{%
\widehat{\Phi}\left(%
\frac{%
2\log\widehat{p}_{k_{u}}%
}{\log(q^r)}%
\right)%
}%
\right)%
\ll\log^{\alpha-1}{(q)}.
\end{multline}
We have altogether
\begin{multline}\label{eq_warwick}
\Eh[q]\left(P_q^2[\Phi;r]^\alpha\right)=\frac{(-2)^\alpha}{\log^\alpha{(q^r)}}\sum_{1\leq j_1,\dotsc,j_\alpha\leq r}(-1)^{\alpha r-(j_1+\dotsc+j_\alpha)} \\
\times\sum_{\substack{k_1,\dotsc,k_\alpha\\\text{distinct}}}\left(\prod_{u=1}^\alpha\left(\frac{\log\widehat{p}_{k_{u}}}{\widehat{p}_{k_{u}}}\widehat{\Phi}\left(\frac{
2\log\widehat{p}_{k_{u}}}{\log(q^r)}\right)\right)\right)\Eh[q]\left(\lambda_f\left(\prod_{u=1}^\alpha\widehat{p}_{k_{u}}^{2j_u}\right)\right) \\
+O\left(\frac{1}{\log{q}}\right)
\end{multline}
since the only element of $P(\alpha,\alpha)$ is the identity function. By lemmas~\ref{lem_delun} and \ref{deltaestimate}, we have
\[%
\Eh[q]\left(\lambda_f\left(\prod_{u=1}^\alpha\widehat{p}_{k_{u}}^{2j_u}\right)\right)%
\ll%
\frac{1}{q}\prod_{u=1}^\alpha \widehat{p}_{k_{u}}^{j_u/2}\log{\widehat{p}_{k_{u}}}
\]
hence the first term in the right-hand side of \eqref{eq_warwick} is bounded by a negative power of $q$ as soon as $\alpha\nu r^2 < 4$.

\subsection{Proof of the third bullet of proposition \ref{moments}} By proposition \ref{explicit}, we have
\begin{equation}\label{eq_EhPqU}
\Eh[q](P_q^1[\Phi;r]^\alpha)%
=%
\frac{(-2)^\alpha}{\log^\alpha{(q^r)}}%
\sum_{\substack{p_1,\dotsc,p_\alpha\in\prem\\ p_1,\dotsc,p_\alpha\nmid q}}%
\left(\prod_{i=1}^\alpha\frac{\log p_i}{\sqrt{p_i}}\widehat{\Phi}\left(\frac{\log p_i}{\log q^r}\right)\right)%
\Eh[q]\left(\prod_{i=1}^\alpha\lambda_f\left(p_i^r\right)\right).%
\end{equation}
Using lemma \ref{primesdistincts}, we rewrite equation \eqref{eq_EhPqU} as
\begin{align}
\Eh[q](P_q^1[\Phi;r]^\alpha)%
&=%
\frac{(-2)^\alpha}{\log^\alpha{(q^r)}}%
\sum_{s=1}^\alpha%
\sum_{\sigma\in P(\alpha,s)}%
\sum_{\substack{i_1,\dotsc,i_s\\\text{distinct}}}%
\left(\prod_{j=1}^\alpha%
\left(\frac{\log\widehat{p}_{i_{\sigma(j)}}}{\sqrt{\widehat{p}_{i_{\sigma(j)}}}}%
\widehat{\Phi}\left(\frac{\log\widehat{p}_{i_{\sigma(j)}}}{\log{(q^r)}}\right)\right)\right)%
\\%
&%
\phantom{%
=\frac{(-2)^\alpha}{\log^\alpha{(q^r)}}%
\sum_{s=1}^\alpha%
\sum_{\sigma\in P(\alpha,s)}%
\sum_{\substack{i_1,\dotsc,i_s\\\text{distinct}}}%
}%
\times%
\Eh[q]\left(\prod_{j=1}^\alpha\lambda_f\left(\widehat{p}_{i_{\sigma(j)}}^r\right)\right)%
\\%
&=%
\frac{(-2)^\alpha}{\log^\alpha{(q^r)}}%
\sum_{s=1}^\alpha%
\sum_{\sigma\in P(\alpha,s)}%
\sum_{\substack{i_1,\dotsc,i_s\\\text{distinct}}}%
\left(%
\prod_{u=1}^s%
\left(%
\frac{%
\log\widehat{p}_{i_u}%
}{%
\sqrt{%
\widehat{p}_{i_u}%
}%
}%
\widehat{\Phi}\left(%
\frac{%
\log\widehat{p}_{i_{u}}%
}{\log q^r}%
\right)%
\right)^{\varpi^{(\sigma)}_u}%
\right)%
\\%
&%
\label{eq_ewrite}
\phantom{%
=
\frac{(-2)^\alpha}{\log^\alpha{(q^r)}}%
\sum_{s=1}^\alpha%
\sum_{\sigma\in P(\alpha,s)}%
\sum_{\substack{i_1,\dotsc,i_s\\\text{distinct}}}%
}
\times%
\Eh[q]\left(%
\prod_{u=1}^s%
\lambda_f\left(\widehat{p}_{i_{u}}^r\right)^{\varpi^{(\sigma)}_u}%
\right).
\end{align}
It follows from \eqref{eq_ams} and \eqref{eq_lintch} that
\[%%
\lambda_f\left(\widehat{p}_{i_{u}}^r\right)^{\varpi^{(\sigma)}_u}%
=%
\sum_{j_u=0}^{r\varpi^{(\sigma)}_u}%
x(\varpi^{(\sigma)}_u,r,j_u)\lambda_f\left(\widehat{p}_{i_{u}}^{j_u}\right).
\]
Since $u\neq v$ implies that $\widehat{p}_{i_{u}}\neq\widehat{p}_{i_{v}}$, equation \eqref{eq_ewrite} becomes
\begin{multline}
\label{eq_drole}
\Eh[q](P_q^1[\Phi;r]^\alpha)%
=%
\frac{(-2)^\alpha}{\log^\alpha{(q^r)}}%
\sum_{s=1}^\alpha%
\sum_{\sigma\in P(\alpha,s)}%
\sum_{\substack{i_1,\dotsc,i_s\\\text{distinct}}}%
\left(%
\prod_{u=1}^s%
\left(%
\frac{%
\log\widehat{p}_{i_u}%
}{%
\sqrt{%
\widehat{p}_{i_u}%
}%
}%
\widehat{\Phi}\left(%
\frac{%
\log\widehat{p}_{i_{u}}%
}{\log{(q^r)}}%
\right)%
\right)^{\varpi^{(\sigma)}_u}%
\right)%
\\%
\times%
\sum_{\substack{j_1,\dotsc,j_s \\
0\leq j_u\leq r\varpi^{(\sigma)}_u}}\left(%
\prod_{u=1}^sx(\varpi^{(\sigma)}_u,r,j_u)%
\right)%
\Eh[q]\left(%
\lambda_f\left(%
\prod_{u=1}^s\widehat{p}_{i_u}^{j_u}%
\right)%
\right).
\end{multline}
%............................CHANGE
Using proposition \ref{iwlusatr} and lemmas \ref{deltaestimate} and \ref{lem_delun}, we get
\[%
\Eh[q]\left(%
\lambda_f\left(%
\prod_{u=1}^s\widehat{p}_{i_u}^{j_u}%
\right)%
\right)%
=%
\prod_{u=1}^s\delta_{j_u,0}+O\left(%
\frac{1}{q}\prod_{u=1}^s\widehat{p}_{i_u}^{j_u/4}\log\widehat{p}_{i_u}
\right)%
\]
hence 
\begin{equation}\label{eq_cesttpte}
\Eh[q](P_q^1[\Phi;r]^\alpha)=\TP+O(\TE)%
\end{equation}
with
\begin{equation}\label{eq_deftp}
\TP\coloneqq%
\frac{(-2)^\alpha}{\log^\alpha{(q^r)}}%
\sum_{s=1}^\alpha%
\sum_{\sigma\in P(\alpha,s)}%
\sum_{\substack{i_1,\dotsc,i_s\\\text{distinct}}}%
\prod_{u=1}^s%
\left(%
\frac{%
\log\widehat{p}_{i_u}%
}{%
\sqrt{%
\widehat{p}_{i_u}%
}%
}%
\widehat{\Phi}\left(%
\frac{%
\log\widehat{p}_{i_{u}}%
}{\log{(q^r)}}%
\right)%
\right)^{\varpi^{(\sigma)}_u}%
x(\varpi^{(\sigma)}_u,r,0)
\end{equation}
and
\begin{equation}\label{eq_defte}
\TE\coloneqq%
\frac{1}{q\log^\alpha{(q^r)}}%
\sum_{s=1}^\alpha%
\sum_{\sigma\in P(\alpha,s)}%
\sum_{\substack{i_1,\dotsc,i_s\\\text{distinct}}}%
\prod_{u=1}^s%
\left(%
\widehat{p}_{i_u}^{(r-2)/4}\log^2\widehat{p}_{i_u}%
\abs{%
\widehat{\Phi}\left(%
\frac{%
\log\widehat{p}_{i_{u}}%
}{\log{(q^r)}}%
\right)%
}
\right)^{\varpi^{(\sigma)}_u}.%
\end{equation}
We have
\begin{equation}\label{eq_majte}
\TE=%
\frac{1}{q\log^\alpha{(q^r)}}%
\left(%
\sum_{\substack{p\in\prem\\ p\nmid q}}p^{(r-2)/4}\log^2p%
\abs{%
\widehat{\Phi}\left(%
\frac{\log p}{\log{(q^r)}}%
\right)%
}%
\right)^\alpha%
\ll%
q^{\alpha r\nu(r+2)/4-1}
\end{equation}
so that, $\TE$ is an error term as soon as
\begin{equation}\label{eq_conddeux}%
\alpha r\nu(r+2)<4.
\end{equation}
We assume from now on that this condition is satisfied. According to \eqref{propriox} (recall that $r\geq 1$), we rewrite \eqref{eq_deftp} as
\begin{equation}\label{eq_tpsimpl}
\TP=%
\frac{(-2)^\alpha}{\log^\alpha{(q^r)}}%
\sum_{s=1}^\alpha%
\sum_{\sigma\in P^{\geq 2}(\alpha,s)}%
\sum_{\substack{i_1,\dotsc,i_s\\\text{distinct}}}%
\prod_{u=1}^s%
\left(%
\frac{%
\log\widehat{p}_{i_u}%
}{%
\sqrt{%
\widehat{p}_{i_u}%
}%
}%
\widehat{\Phi}\left(%
\frac{%
\log\widehat{p}_{i_{u}}%
}{\log{(q^r)}}%
\right)%
\right)^{\varpi^{(\sigma)}_u}%
x(\varpi^{(\sigma)}_u,r,0)
\end{equation}
where
\[%
P^{\geq 2}(\alpha,s)\coloneqq%
\left\{%
\sigma\in P(\alpha,s) \colon %
\forall u\in\{1,\dotsc,s\}, % 
\varpi^{(\sigma)}_u\geq 2%
\right\}.
\]
Moreover, if for at least one $\sigma$ and at least one $u$ (say $u_0$) we have $\varpi^{(\sigma)}_u\geq 3$, then
\begin{multline}
\sum_{\substack{i_1,\dotsc,i_s\\\text{distinct}}}%
\prod_{u=1}^s%
\left(%
\frac{%
\log\widehat{p}_{i_u}%
}{%
\sqrt{%
\widehat{p}_{i_u}%
}%
}%
\widehat{\Phi}\left(%
\frac{%
\log\widehat{p}_{i_{u}}%
}{\log{(q^r)}}%
\right)%
\right)^{\varpi^{(\sigma)}_u}%
x(\varpi^{(\sigma)}_u,r,0)%
\\%
\ll%
\left(%
\sum_{\substack{p\in\prem\\ p\leq q^{r\nu}}}\frac{\log^3{(p)}}{p^{3/2}}%
\right)%
\prod_{\substack{u=1\\ u\neq u_0}}^s
\left(%
\sum_{\substack{p_u\in\prem\\ p_u\leq q^{r\nu}}}\frac{\log^2{(p_u)}}{p_u}%
\right)%
\\%
\ll\label{eq_setm}%
(\log q)^{2s-2}.
\end{multline}
But, from \eqref{eq_trivcond}, we deduce
\[%
2s\leq\sum_{j=1}^s\varpi^{(\sigma)}_j=\alpha%
\]
hence $(\log q)^{2s-2}\ll (\log q)^{\alpha-2}$. Reinserting this in \eqref{eq_setm} and the result in \eqref{eq_tpsimpl}, we obtain %
\begin{multline}\label{eq_tpplussimpl}
\TP=%
\frac{(-2)^\alpha}{\log^\alpha{(q^r)}}%
\sum_{s=1}^\alpha%
\sum_{\sigma\in P^{2}(\alpha,s)}%
\sum_{\substack{i_1,\dotsc,i_s\\\text{distinct}}}%
\prod_{u=1}^s%
\left(%
\frac{%
\log\widehat{p}_{i_u}%
}{%
\sqrt{%
\widehat{p}_{i_u}%
}%
}%
\widehat{\Phi}\left(%
\frac{%
\log\widehat{p}_{i_{u}}%
}{\log q^r}%
\right)%
\right)^{\varpi^{(\sigma)}_u}%
x(\varpi^{(\sigma)}_u,r,0) \\
+O\left(\frac{1}{\log^2{(q)}}\right)
\end{multline}
where
\[%
P^{2}(\alpha,s)\coloneqq%
\left\{%
\sigma\in P(\alpha,s) \colon %
\forall u\in\{1,\dotsc,s\}, % 
\varpi^{(\sigma)}_u= 2%
\right\}.
\]
From \eqref{eq_tpplussimpl}, \eqref{eq_majte} and \eqref{eq_cesttpte}, we deduce
%.....................ENDCHANGE
\begin{multline}
\Eh[q](P_q^1[\Phi;r]^\alpha)%
=%
\frac{(-2)^\alpha}{\log^\alpha{(q^r)}}%
\sum_{s=1}^\alpha%
\sum_{\sigma\in P^{2}(\alpha,s)}%
\sum_{\substack{i_1,\dotsc,i_s\\\text{distinct}}}%
\prod_{u=1}^s%
\frac{%
\log^2{(\widehat{p}_{i_u})}%
}{%%
\widehat{p}_{i_u}%
}%
\widehat{\Phi}^2\left(%
\frac{%
\log\widehat{p}_{i_{u}}%
}{\log{(q^r)}}%
\right) \\
+O\left(\frac{1}{\log^2{(q)}}\right)
\end{multline}
since $x(2,r,0)=1$ according to \eqref{propriox}. Note in particular that, according to \eqref{eq_trivcond} the previous sum is zero if $\alpha$ is odd. Thus, we can assume now that $\alpha$ is even and get
\begin{multline}\label{eq_ouf}
\Eh[q](P_q^1[\Phi;r]^\alpha)%
=%
\frac{(-2)^\alpha}{\log^\alpha{(q^r)}}%
\sum_{\sigma\in P^{2}(\alpha,\alpha/2)}%
\sum_{\substack{i_1,\dotsc,i_{\alpha/2}\\\text{distinct}}}%
\prod_{u=1}^{\alpha/2}%
\frac{%
\log^2{(\widehat{p}_{i_u})}%
}{%%
\widehat{p}_{i_u}%
}%
\widehat{\Phi}^2\left(%
\frac{%
\log\widehat{p}_{i_{u}}%
}{\log{(q^r)}}%
\right) \\
+O\left(\frac{1}{\log^2{(q)}}\right).
\end{multline}
However, summing over all the possible $(i_1,\dotsc,i_{\alpha/2})$ instead of the one with distinct indices reintroduces
convergent sums that enter the error term because of the $1/\log^\alpha{(q^r)}$ factor. It follows that \eqref{eq_ouf} becomes:
\begin{multline}\label{eq_ooo}
\Eh[q](P_q^1[\Phi;r]^\alpha)%
=%
\left[%
\frac{4}{\log^2{(q^r)}}%
\sum_{p\in\prem}%
\frac{\log^2{(p)}}{p}%
\widehat{\Phi}^2\left(%
\frac{\log p}{\log{(q^r)}}%
\right)%
\right]^{\alpha/2}%
\#P^{2}(\alpha,\alpha/2) \\
+O\left(\frac{1}{\log^2{(q)}}\right).
\end{multline}
Taking $m=2$ (we already proved that the second moment is finite, see section~\ref{sec_twoandvar}) %
and reinserting the result in \eqref{eq_ooo} implies that
\[%%
\Eh[q](P_q^1[\Phi;r]^\alpha)%
=%
\Eh[q](P_q^1[\Phi;r]^2)%
\#P^{2}(\alpha,\alpha/2)+O\left(\frac{1}{\log^2{(q)}}\right).
\]
We conclude by computing
\[%
\#P^{2}(\alpha,\alpha/2)=\frac{\alpha!}{2^{\alpha/2}\left(\frac{\alpha}{2}\right)!}.
\]
(see \cite[Example 5.2.6 and Exercise 5.43]{Stan00}).
\subsection{Proof of the second bullet of proposition \ref{moments}} %
We mix the two techniques which have been used to prove the first and third %
bullets of proposition \ref{moments}. We get following the same lines and %
thanks to lemma~\ref{primesdistincts}
\begin{multline}
\Eh[q]\left(P_q^1[\Phi;r]^{m-\ell}P_q^2[\Phi;r]^\alpha\right)=%
\frac{(-2)^{\alpha+m-\ell}}{\log^{\alpha+m-\ell}{(q^r)}}%
\sum_{1\leq j_1,\dotsc,j_\alpha\leq r}(-1)^{\alpha r-(j_1+\dotsc+j_\alpha)}%
\sum_{s=1}^{\alpha+m-\ell} \\
\times%
\sum_{\sigma\in%
  P(\alpha+m-\ell,s)}\sum_{\substack{i_1,\dotsc,i_s\\\text{distinct}}}%
\prod_{u=1}^{s}%
\left(%
\frac{\log^{\varpi^{(\sigma,1)}_{u}+\varpi^{(\sigma,2)}_{u}}%
{\left(\widehat{p}_{i_u}\right)}}{\widehat{p}_{i_{u}}^{\varpi^{(\sigma,1)}_{u}/2+\varpi^{(\sigma,2)}_{u}}}%
\widehat{\Phi}\left(\frac{\log\widehat{p}_{i_{u}}}{\log(q^r)}\right)^{\varpi^{(\sigma,1)}_{u}}%
\widehat{\Phi}\left(\frac{2\log\widehat{p}_{i_{u}}}{\log(q^r)}\right)^{\varpi^{(\sigma,2)}_{u}}%
\right) \\
\times %
\Eh[q]\left(\prod_{u=1}^s\left(\lambda_f\left(\widehat{p}_{i_{u}}^{r}\right)^{\varpi^{(\sigma,1)}_{u}}%
\prod_{j=1}^r\lambda_f\left(\widehat{p}_{i_{u}}^{2j}\right)^{\varpi^{(\sigma,2)}_{u,j}}\right)\right)
\end{multline}
where
\begin{align*}
\varpi^{(\sigma,1)}_{u} & \coloneqq \#\left\{i\in\left\{1,\dotsc,m-\ell\right\},\; \sigma(i)=u\right\}, \\
\varpi^{(\sigma,2)}_{u} & \coloneqq \#\left\{i\in\left\{1,\dotsc,\alpha\right\},\; \sigma(m-\ell+i)=u\right\}, \\
\varpi^{(\sigma,2)}_{u,j} & \coloneqq %
\#\left\{i\in\left\{1,\dotsc,\alpha\right\},\; \sigma(m-\ell+i)=u \text{ and } j_{i}=j\right\}
%\#\left\{i\in\left\{1,\dotsc,\alpha\right\},\; \sigma(m-\ell+i)=u \text{ and } j_{m-\ell+i}=j\right\}
\end{align*}
for any $1\leq u\leq s$, any $1\leq j\leq r$ and any $\sigma\in P(\alpha+m-\ell,s)$. Note that these numbers satisfy
\begin{equation}
\label{prop1}
\sum_{u=1}^s\left(\varpi^{(\sigma,1)}_{u}+\varpi^{(\sigma,2)}_{u}\right)=m-\ell+\alpha
\end{equation}
and
\begin{equation}
\label{prop2}
\sum_{j=1}^r\varpi^{(\sigma,2)}_{u,j}=\varpi^{(\sigma,2)}_{u}
\end{equation}
for any $1\leq u\leq r$ and any $\sigma\in P(\alpha+m-\ell,s)$ by definition. They also satisfy
\begin{equation}
\label{prop3}
\forall\sigma\in P(\alpha+m-\ell,s), \forall u\in\left\{1,\dotsc,s\right\}, \quad  \varpi^{(\sigma,1)}_{u}+\varpi^{(\sigma,2)}_{u}\geq 1
\end{equation}
since any $\sigma\in P(\alpha+m-\ell,s)$ is surjective and
\begin{equation}
\label{prop4}
\forall\sigma\in P(\alpha+m-\ell,s), \forall i\in\left\{1,2\right\},\exists u_{i,\sigma}\in\left\{1,\dotsc,s\right\}, \quad  \varpi^{(\sigma,i)}_{u_{i,\sigma}}\geq 1
\end{equation}
since $\alpha\geq 1$ and $m-\ell\geq 1$. The strategy is to estimate individually each term of the $\sigma$-sum. Thus, we fix some integers $j_1,\dotsc,j_\alpha$ in $\left\{1,\dotsc,r\right\}$, some integer $s$ in $\left\{1,\dotsc,r\right\}$ and some application $\sigma$ in $P(\alpha+m-\ell,s)$.\newline
\noindent{\textit{\underline{First case}:} $\quad\mathit{\forall u\in\left\{1,\dotsc,s\right\}, \;\varpi^{(\sigma,1)}_{u}/2+\varpi^{(\sigma,2)}_{u}\leq 1.}$}\newline
Let us remark that if $\varpi^{(\sigma,2)}_{u}=1$ for some $1\leq u\leq s$ then there exists a unique $1\leq j_{i_u}\leq r$ depending on $\sigma$ such that $\varpi^{(\sigma,2)}_{u,j_{i_u}}=1$ and $\varpi^{(\sigma,2)}_{u,j}=0$ for any $1\leq j\neq j_{i_u}\leq r$ according to \eqref{prop2}. Thus,
\begin{multline*}
\prod_{u=1}^s\left(\lambda_f\left(\widehat{p}_{i_{u}}^{r}\right)^{\varpi^{(\sigma,1)}_{u}}\prod_{j=1}^r\left(\lambda_f\left(\widehat{p}_{i_{u}}^{2j}\right)^{\varpi^{(\sigma,2)}_{u,j}}\right)\right)=\lambda_f\left(\prod_{\substack{1\leq u\leq s\\
\left(\varpi^{(\sigma,1)}_{u},\varpi^{(\sigma,2)}_{u}\right)=(2,0)}}\widehat{p}_{i_{u}}^{r\varpi^{(\sigma,1)}_{u}/2}\right) \\
\times\lambda_f\left(\prod_{\substack{1\leq u\leq s\\
\left(\varpi^{(\sigma,1)}_{u},\varpi^{(\sigma,2)}_{u}\right)=(2,0)}}\widehat{p}_{i_{u}}^{r\varpi^{(\sigma,1)}_{u}/2}\prod_{\substack{1\leq u\leq s\\
\left(\varpi^{(\sigma,1)}_{u},\varpi^{(\sigma,2)}_{u}\right)=(1,0)}}\widehat{p}_{i_{u}}^{r\varpi^{(\sigma,1)}_{u}}\prod_{\substack{1\leq u\leq s\\
\left(\varpi^{(\sigma,1)}_{u},\varpi^{(\sigma,2)}_{u}\right)=(0,1)}}\widehat{p}_{i_{u}}^{2j_{i_u}\varpi^{(\sigma,2)}_{u,j_{i_u}}}\right)
\end{multline*}
where the two integers appearing in the right-hand side of the previous equality are different according to \eqref{prop4}. %
Consequently, proposition \ref{iwlusatr} and lemmas \ref{deltaestimate} and \ref{lem_delun} enable us to assert that %
\begin{multline*}
\Eh[q]\left(\prod_{u=1}^s\left(\lambda_f\left(\widehat{p}_{i_{u}}^{r}\right)^{\varpi^{(\sigma,1)}_{u}}\prod_{j=1}^r\left(\lambda_f\left(\widehat{p}_{i_{u}}^{2j}\right)^{\varpi^{(\sigma,2)}_{u,j}}\right)\right)\right) 
\ll%
\frac{1}{q}%
\prod_{\substack{1\leq u\leq s\\ \left(\varpi^{(\sigma,1)}_{u},\varpi^{(\sigma,2)}_{u}\right)=(2,0)}}%
                  \frac{\log{\widehat{p}_{i_{u}}}}{\widehat{p}_{i_{u}}^{-r\varpi^{(\sigma,1)}_{u}/4}}%
\\
\times%
\prod_{\substack{1\leq u\leq s\\ \left(\varpi^{(\sigma,1)}_{u},\varpi^{(\sigma,2)}_{u}\right)=(1,0)}}%
                  \frac{\log{\widehat{p}_{i_{u}}}}{\widehat{p}_{i_{u}}^{-r\varpi^{(\sigma,1)}_{u}/4}}%
\prod_{\substack{1\leq u\leq s\\ \left(\varpi^{(\sigma,1)}_{u},\varpi^{(\sigma,2)}_{u}\right)=(0,1)}}%
                  \frac{\log{\widehat{p}_{i_{u}}}}{\widehat{p}_{i_{u}}^{-r\varpi^{(\sigma,2)}_{u}/2}}.%
\end{multline*}
Note that, in this first case, the right hand term is %
\[%
\frac{1}{q}%
\prod_{u=1}^s%
\frac{%
  \log{\widehat{p}_{i_{u}}}%
     }%
     {%
  \widehat{p}_{i_{u}}^{-r(\varpi^{(\sigma,1)}_{u}/4+\varpi^{(\sigma,2)}_{u}/2)}%
     }
\]
hence the contribution of these $\sigma$'s to
$\Eh[q]\left(P_q^1[\Phi;r]^{m-\ell}P_q^2[\Phi;r]^\alpha\right)$ is bounded by %
\[%
\frac{q^\epsilon}{q}
\left(%
 \sum_{p\leq q^{\nu r}}\frac{1}{p^{1/2-r/4}}%
\right)^{m-\ell}%
\left(%
 \sum_{p\leq q^{\nu r/2}}\frac{1}{p^{1-r/2}}%
\right)^{\alpha}%
\ll%
q^{%
\nu r/4[(m-\ell)(r+2)+\alpha r]-1+\epsilon%
}.
\]
This is an admissible error term as long as %
$\nu r/4[(m-\ell)(r+2)+\alpha r]<1$.

\noindent{\textit{\underline{Second case}:} $\quad\mathit{\exists
    u_\sigma\in\left\{1,\dotsc,s\right\},
    \;\varpi^{(\sigma,1)}_{u_\sigma}/2+\varpi^{(\sigma,2)}_{u_\sigma}>1.}$}

According to \eqref{eq_ams} and \eqref{eq_lintch}, if $1\leq u\leq s$ and $1\leq j\leq r$ then
\begin{equation*}
\lambda_f\left(\widehat{p}_{i_{u}}^r\right)^{\varpi^{(\sigma,1)}_u}=\sum_{k_{u,1}=0}^{r\varpi^{(\sigma,1)}_u}x(\varpi^{(\sigma,1)}_u,r,k_{u,1})\lambda_f\left(\widehat{p}_{i_{u}}^{k_{u,1}}\right)
\end{equation*}
and
\begin{equation*}
\lambda_f\left(\widehat{p}_{i_{u}}^{2j}\right)^{\varpi^{(\sigma,2)}_{u,j}}=%
\sum_{k_{u,j,2}=0}^{j\varpi^{(\sigma,2)}_{u,j}}%
x(\varpi^{(\sigma,2)}_{u,j},2j,2k_{u,j,2})%
\lambda_f\left(\widehat{p}_{i_{u}}^{2k_{u,j,2}}\right)
\end{equation*}
since $x(\varpi^{(\sigma,2)}_{u,j},2j,k_{u,j,2})=0$ if $k_{u,j,2}$ is odd (see \eqref{propriox}). Then, one may remark that
\begin{equation*}
\prod_{1\leq j\leq r}\lambda_f\left(\widehat{p}_{i_{u}}^{2k_{u,j,2}}\right)=\sum_{\ell_u=0}^{K_u}y_{\ell_u}\lambda_f\left(\widehat{p}_{i_{u}}^{2\ell_u}\right)
\end{equation*}
for some integers $y_{\ell_u}$ and where $K_u\coloneqq\sum_{1\leq j\leq r}k_{u,j,2}$ for any $1\leq u\leq s$. All these facts lead to
\begin{multline}
\Eh[q]\left(P_q^1[\Phi;r]^{m-\ell}P_q^2[\Phi;r]^\alpha\right)=\frac{(-2)^{\alpha+m-\ell}(-1)^{\alpha r}}{\log^{\alpha+m-\ell}{(q^r)}}\sum_{1\leq j_1,\dotsc,j_\alpha\leq r}(-1)^{j_1+\dotsc+j_\alpha}\sum_{s=1}^{\alpha+m-\ell} \\
\times\sum_{\sigma\in P(\alpha+m-\ell,s)}\sum_{\substack{i_1,\dotsc,i_s\\\text{distinct}}}\prod_{u=1}^{s}\left(\frac{\log^{\varpi^{(\sigma,1)}_{u}+\varpi^{(\sigma,2)}_{u}}{\left(\widehat{p}_{i_u}\right)}}{\widehat{p}_{i_{u}}^{\varpi^{(\sigma,1)}_{u}\left/2\right.+\varpi^{(\sigma,2)}_{u}}}\widehat{\Phi}\left(\frac{\log\widehat{p}_{i_{u}}}{\log(q^r)}\right)^{\varpi^{(\sigma,1)}_{u}}\widehat{\Phi}\left(\frac{
2\log\widehat{p}_{i_{u}}}{\log(q^r)}\right)^{\varpi^{(\sigma,2)}_{u}}\right) \\
\times\sum_{\substack{0\leq k_{1,1}\leq r\varpi^{(\sigma,1)}_1 \\
\vdots \\
0\leq k_{s,1}\leq r\varpi^{(\sigma,1)}_s}}\sum_{\substack{0\leq k_{1,1,2}\leq \varpi^{(\sigma,2)}_{1,1} \\
\vdots \\
0\leq k_{s,1,2}\leq\varpi^{(\sigma,2)}_{s,1}}}\dotsc\sum_{\substack{0\leq k_{1,r,2}\leq r\varpi^{(\sigma,2)}_{1,r} \\
\vdots \\
0\leq k_{s,r,2}\leq r\varpi^{(\sigma,2)}_{s,r}}}\sum_{\substack{0\leq\ell_1\leq K_1 \\
\vdots \\
0\leq\ell_s\leq K_s}} \\
\times\prod_{u=1}^s\left(x\left(\varpi^{(\sigma,1)}_u,r,k_{u,1}\right)y_{\ell_u}\prod_{j=1}^r\left(x\left(\varpi^{(\sigma,2)}_{u,j},2j,2k_{u,j,2}\right)\right)\right) \\
\times\Eh[q]\left(\lambda_f\left(\prod_{u=1}^s\widehat{p}_{i_{u}}^{k_{u,1}}\right)\lambda_f\left(\prod_{u=1}^s\widehat{p}_{i_{u}}^{2\ell_u}\right)\right).
\end{multline}
Proposition \ref{iwlusatr} and lemmas \ref{deltaestimate} and \ref{lem_delun} enable us to assert that
\begin{equation*}
\Eh[q]\left(\lambda_f\left(\prod_{u=1}^s\widehat{p}_{i_{u}}^{k_{u,1}}\right)\lambda_f\left(\prod_{u=1}^s\widehat{p}_{i_{u}}^{2\ell_u}\right)\right)=\prod_{u=1}^s\delta_{k_{u,1},2\ell_u} \\
+O\left(\frac{1}{q}\prod_{u=1}^s\widehat{p}_{i_{u}}^{k_{u,1}/4+\ell_u/2}\log{\widehat{p}_{i_{u}}}\right)
\end{equation*}
and we can write
\begin{equation} 
\Eh[q]\left(P_q^1[\Phi;r]^{m-\ell}P_q^2[\Phi;r]^\alpha\right)=\TP+O(\TE)%
\end{equation}
with
\begin{multline}\label{endddd}
\TP\coloneqq\frac{(-2)^{\alpha+m-\ell}(-1)^{\alpha r}}{\log^{\alpha+m-\ell}{(q^r)}}\sum_{1\leq j_1,\dotsc,j_\alpha\leq r}(-1)^{j_1+\dotsc+j_\alpha}\sum_{s=1}^{\alpha+m-\ell} \\
\times\sum_{\sigma\in P(\alpha+m-\ell,s)}\sum_{\substack{i_1,\dotsc,i_s\\\text{distinct}}}\prod_{u=1}^{s}\left(\frac{\log^{\varpi^{(\sigma,1)}_{u}+\varpi^{(\sigma,2)}_{u}}{\left(\widehat{p}_{i_u}\right)}}{\widehat{p}_{i_{u}}^{\varpi^{(\sigma,1)}_{u}\left/2\right.+\varpi^{(\sigma,2)}_{u}}}\widehat{\Phi}\left(\frac{\log\widehat{p}_{i_{u}}}{\log(q^r)}\right)^{\varpi^{(\sigma,1)}_{u}}\widehat{\Phi}\left(\frac{
2\log\widehat{p}_{i_{u}}}{\log(q^r)}\right)^{\varpi^{(\sigma,2)}_{u}}\right) \\
\times\sum_{\substack{0\leq k_{1,1,2}\leq \varpi^{(\sigma,2)}_{1,1} \\
\vdots \\
0\leq k_{s,1,2}\leq\varpi^{(\sigma,2)}_{s,1}}}\dotsc\sum_{\substack{0\leq k_{1,r,2}\leq r\varpi^{(\sigma,2)}_{1,r} \\
\vdots \\
0\leq k_{s,r,2}\leq r\varpi^{(\sigma,2)}_{s,r}}}\sum_{\substack{0\leq\ell_1\leq r\min{\left(\varpi^{(\sigma,1)}_{1}/2,\varpi^{(\sigma,2)}_{1}\right)} \\
\vdots \\
0\leq\ell_s\leq r\min{\left(\varpi^{(\sigma,1)}_{s}/2,\varpi^{(\sigma,2)}_{s}\right)}}} \\
\times\prod_{u=1}^{s}\left(x\left(\varpi^{(\sigma,1)}_u,r,2\ell_u\right)y_{\ell_u}\prod_{j=1}^r\left(x\left(\varpi^{(\sigma,2)}_{u,j},2j,2k_{u,j,2}\right)\right)\right)
\end{multline}
and
\begin{multline}
\TE\coloneqq%
\frac{1}{q\log^{\alpha+m-\ell}{(q^r)}}%
\\%
\times%
\sum_{s=1}^{\alpha+m-\ell}%
\sum_{\sigma\in P(\alpha+m-\ell,s)}%
\sum_{\substack{i_1,\dotsc,i_s\\\text{distinct}}}
\prod_{u=1}^{s}\log^{\varpi^{(\sigma,1)}_{u}+\varpi^{(\sigma,2)}_{u}+1}{(\widehat{p}_{i_u})}
\widehat{p}_{i_u}^{\left(r/2-1\right)\left(\varpi^{(\sigma,1)}_{u}/2+\varpi^{(\sigma,2)}_{u}\right)}%
\\%
\times%
\left\vert\widehat{\Phi}\left(\frac{\log\widehat{p}_{i_{u}}}{\log(q^r)}\right)\right\vert^{\varpi^{(\sigma,1)}_{u}}%
\left\vert\widehat{\Phi}\left(\frac{2\log\widehat{p}_{i_{u}}}{\log(q^r)}\right)\right\vert^{\varpi^{(\sigma,2)}_{u}}
\end{multline}
which is bounded by $O_\epsilon\left(q^{(\alpha+m-\ell)\nu r^2/4-1+\varepsilon}\right)$ for any $\varepsilon>0$ and is an admissible error term if $(\alpha+m-\ell)\nu<4/r^2$. Estimating $\TP$ is possible since we can assume that $\sigma$ satisfies the following additional property. If $\varpi_u^{(\sigma,2)}=0$ for some $1\leq u\leq s$ then $\varpi_u^{(\sigma,1)}>1$. Let us assume on the contrary that $\varpi_u^{(\sigma,1)}\leq 1$ which entails $\varpi_u^{(\sigma,1)}=1$ according to \eqref{prop3}. Then,
\begin{equation*}
x\left(\varpi^{(\sigma,1)}_u,r,2\ell_u\right)=x\left(1,r,0\right)=0
\end{equation*} 
since $\ell_u=0$ and according to \eqref{propriox}. Thus, the contribution of the $\sigma$'s which do not satisfy this last property vanishes. As a consequence, the sum over the distinct $i_1,\dotsc,i_s$ is bounded by
\begin{multline*}
\sum_{\substack{i_1,\dotsc,i_s\\\text{distinct}}}\prod_{\substack{1\leq u\leq s \\
\left(\varpi_u^{(\sigma,1)},\varpi_u^{(\sigma,2)}\right)=(2,0)}}\left(\frac{\log^{2}{\left(\widehat{p}_{i_u}\right)}}{\widehat{p}_{i_{u}}}\left\vert\widehat{\Phi}\left(\frac{\log\widehat{p}_{i_{u}}}{\log(q^r)}\right)\right\vert^{2}\right) \\
\times\prod_{\substack{1\leq u\leq s \\
\left(\varpi_u^{(\sigma,1)},\varpi_u^{(\sigma,2)}\right)=(0,1)}}\left(\frac{\log{\left(\widehat{p}_{i_u}\right)}}{\widehat{p}_{i_{u}}}\left\vert\widehat{\Phi}\left(\frac{2\log\widehat{p}_{i_{u}}}{\log(q^r)}\right)\right\vert\right) \\
\times\prod_{\substack{1\leq u\leq s \\
\varpi_u^{(\sigma,1)}/2+\varpi_u^{(\sigma,2)}>1}}\left(\frac{\log^{\varpi^{(\sigma,1)}_{u}+\varpi^{(\sigma,2)}_{u}}{\left(\widehat{p}_{i_u}\right)}}{\widehat{p}_{i_{u}}^{\varpi^{(\sigma,1)}_{u}\left/2\right.+\varpi^{(\sigma,2)}_{u}}}\left\vert\widehat{\Phi}\left(\frac{\log\widehat{p}_{i_{u}}}{\log(q^r)}\right)\right\vert^{\varpi^{(\sigma,1)}_{u}}\left\vert\widehat{\Phi}\left(\frac{
2\log\widehat{p}_{i_{u}}}{\log(q^r)}\right)\right\vert^{\varpi^{(\sigma,2)}_{u}}\right)
\end{multline*}
which is itself bounded by $O\left(\log^{A_\sigma}{(q)}\right)$ where the exponent is given by
\begin{multline*}
A_\sigma\coloneqq 2\#\left\{1\leq u\leq s, \varpi_u^{(\sigma,2)}=0 \text{ and }  \varpi_u^{(\sigma,1)}/2+\varpi_u^{(\sigma,2)}\leq 1\right\} \\
+\#\left\{1\leq u\leq s, \varpi_u^{(\sigma,2)}=1 \text{ and }  \varpi_u^{(\sigma,1)}/2+\varpi_u^{(\sigma,2)}\leq 1\right\}<m-\ell+\alpha.
\end{multline*}
The last inequality follows from (see \eqref{prop1} and the additional property of $\sigma$)
\begin{equation*}
m-\ell+\alpha=A_\sigma+\sum_{\substack{1\leq u\leq s \\
\varpi_u^{(\sigma,1)}/2+\varpi_u^{(\sigma,2)}>1}}\left(\varpi_u^{(\sigma,1)}+\varpi_u^{(\sigma,2)}\right).
\end{equation*}
Thus, the contribution of the $\TP$ term of these $\sigma$'s to $\Eh[q]\left(P_q^1[\Phi;r]^{m-\ell}P_q^2[\Phi;r]^\alpha\right)$ is bounded by $O\left(\log^{-1}{(q)}\right)$.

%..................

\appendix

%\section{Chebyshev}

\section{Analytic and arithmetic toolbox}
\label{klooster}
\subsection{On smooth dyadic partitions of unity}

\label{unity}

Let $\psi\colon\R_+\to\R$ be any smooth function satisfying
\begin{equation*}
\psi(x)=\begin{cases}
0 & \text{if } 0\leq x\leq 1, \\
1 & \text{if } x>\sqrt{2}
\end{cases}
\end{equation*}
and $x^j\psi^{(j)}(x)\ll_j 1$ for any real number $x\geq 0$ and any integer $j\geq 0$. If $\rho\colon\R_+\to\R$ is the function defined by
\begin{equation*}
\rho(x)\coloneqq\begin{cases}
\psi(x) & \text{if } 0\leq x\leq\sqrt{2}, \\
1-\psi\left(\frac{x}{\sqrt{2}}\right) & \text{otherwise}
\end{cases}
\end{equation*}
then $\rho$ is a smooth function compactly supported in $[1,2]$ satisfying
\begin{equation*}
x^j\rho^{(j)}(x)\ll_j 1\quad\text{ and }\quad\sum_{a\in\mathbb{Z}}\rho\left(\frac{x}{\sqrt{2}^a}\right)=1
\end{equation*}
for any real number $x\geq 0$ and any integer $j\geq 0$. 
 
%\ifpdf
% \begin{figure}[htb]%
%   \centering
%   \subfigure[Graph of $\psi$]{%
%    \includegraphics{GraphePsi.pdf}%
%   }\hfill
%  \subfigure[Graph of $\rho$]{%
%   \includegraphics{GrapheRho.pdf}%
%   }
% \end{figure}
%\else
 \begin{figure}[htb]%
   \centering
   \subfigure[Graph of $\psi$]{%
    \includegraphics{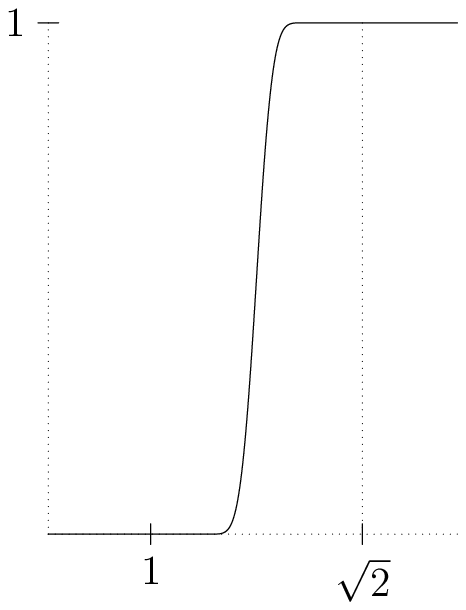}%
   }\hfill
   \subfigure[Graph of $\rho$]{%
   \includegraphics{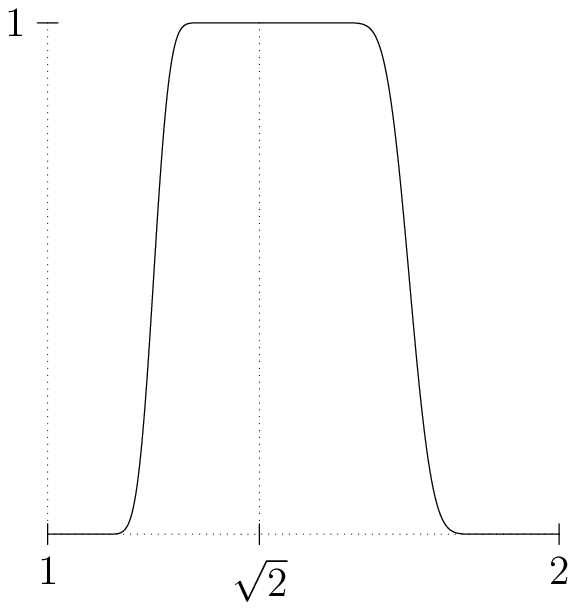}%
   }
 \end{figure}
%\fi

If $F\colon\R_+^n\to\R$ is a function of $n\geq 1$ real variables then we can decompose it in
\begin{equation*}
F=\sum_{a_1\in\mathbb{Z}}\ldots\sum_{a_n\in\mathbb{Z}}F_{A_1,\cdots,A_n}
\end{equation*}
where $A_i\coloneqq\sqrt{2}^{a_i}$ and
\begin{equation*}
F_{A_1,\cdots,A_n}(x_1,\cdots,x_n)\coloneqq\prod_{i=1}^n\rho_{A_i}(x_i)F(x_1,\cdots,x_n)
\end{equation*}
with $\rho_{A_i}(x_i)\coloneqq\rho\left(x_i\left/A_i\right.\right)$ is a smooth function %
compactly supported in $[A_i,2A_i]$ satisfying $x_i^j\rho_{A_i}^{(j)}(x_i)\ll_j 1$ for any real %
number $x_i\geq 0$ and any integer $j\geq 0$. Let us introduce the following notation for summation over powers of $\sqrt{2}$ :
\[%
\sumsh_{A\leq M\leq B}f(M)\coloneqq%
\sum_{%
\substack{n\in\N\\ A\leq 2^{n/2}\leq B}%
}%
f\left(%
2^{n/2}%
\right).%
\]

We will use such smooth dyadic partitions of unity several %
times in this paper and we will also need these natural estimates in such contexts
\begin{equation}
\label{dyadic1}
\sumsh_{M\leq M_1}
M^\alpha%
\ll %
M_1^\alpha%
\end{equation}
for any $\alpha, M_1>0$ and
\begin{equation}
\label{dyadic2}
\sumsh_{M\geq M_0}%
\frac{1}{M^{\alpha}}%
\ll%
\frac{1}{M_0^\alpha}%
\end{equation}
for any $\alpha, M_0>0$.

\subsection{On Bessel functions}
\label{Bessel}

The Bessel function of first kind and order a integer $\kappa\geq 1$ is defined by
\[%
\forall z\in\C, \quad J_\kappa(z)\coloneqq\sum_{n\geq 0}\frac{(-1)^n}{n!(\kappa+n)!}\left(\frac{z}{2}\right)^{\kappa+2n}.
\]

It satisfies the following estimate (founded in \cite[Lemma C.2]{MR1915038}), valid %
for any real number $x$, any integer $j\geq 0$ and any integer $\kappa\geq 1$: %
\begin{equation}
\label{bessel}
\left(\frac{x}{1+x}\right)^jJ_\kappa^{(j)}(x)\ll_{j,\kappa}\frac{1}{\left(1+x\right)^{\frac{1}{2}}}\left(\frac{x}{1+x}\right)^\kappa
\end{equation}
for any real number $x$, any integer $j\geq 0$ and any integer $\kappa\geq 1$. %
The following useful lemma follows immediately. %
\begin{lemma}\label{lem_picard}%
Let $X>0$ and $\kappa\geq 1$, then %
\[%
\sum_{d>0}\frac{\tau(d)}{\sqrt{d}}%
\abs{J_{\kappa}\left(\frac{X}{d}\right)}%
\ll%
\begin{cases}%
X^{1/2}\log X & \text{if $\;X>1$,}\\
X^{\kappa} & \text{if $\;0<X\leq 1$.}
\end{cases}
\]%
\end{lemma}%
\subsection{Basic facts on Kloosterman sums}\label{Kloos}%
%%
%.................
%\subsubsection{Basic facts}%
For any integer $m,n, c\geq 1$, the Kloosterman sum is defined by %
\[%
S(m,n;c)\coloneqq%
\sum_{\substack{x\mod (c) \\ (x,c)=1}}%
e\left(\frac{mx+n\overline{x}}{c}\right)%
\]
where $\overline{x}$ stands for the inverse of $x$ modulo $c$. We recall %
some basic facts on these sums. The Chinese remainder theorem implies the %
following multiplicativity relation %
\begin{equation}\label{eq_crt}%
S(m,n;qr)%
=%
S(m\overline{q}^2,n;r)%
S(m\overline{r}^2,n;q)%
\end{equation}
valid as soon as $(q,r)=1$. Here, $\overline{q}$ (resp. $\overline{r}$) %
is the inverse of $q$ (resp. $r$) modulo $r$ (resp. $q$). If $p$ and $q$ are two prime numbers, $\gamma\geq 1$ and $r\geq 1$ then, from %
\eqref{eq_crt} and \cite[(2.312)]{55.0703.02} we obtain %
\begin{equation}\label{eq_klzero}%
S\left(p^{\gamma}q,1;qr\right)%
=%%
\begin{cases}%
-S\left(p^{\gamma}\overline{q},1;r\right) & \text{if $(q,r)=1$, }\\%
0 & \text{otherwise.}%
\end{cases}%
\end{equation}%
The Weil-Estermann inequality \cite{Es} is %
\begin{equation}\label{weil}%
\abs{S(m,n;c)}\leq\sqrt{(m,n,c)}\tau(c)\sqrt{c}. %
\end{equation}%

%......................................
\providecommand{\bysame}{\leavevmode\hbox to3em{\hrulefill}\thinspace}
\providecommand{\MR}{\relax\ifhmode\unskip\space\fi MR }
\providecommand{\MRhref}[2]{%
  \href{http://www.ams.org/mathscinet-getitem?mr=#1}{#2}
}
\providecommand{\href}[2]{#2}

%.......................................
%\bibliographystyle{amsplain}
%\bibliography{riro}

\begin{thebibliography}{10}

\bibitem{AtLe}
A.~O.~L. Atkin and J.~Lehner, \emph{Hecke operators on {$\Gamma \sb{0}(m)$}},
  Math. Ann. \textbf{185} (1970), 134--160. \MR{MR0268123 (42 \#3022)}

\bibitem{BlHaMi}
V.~Blomer, G.~Harcos, and P.~Michel, \emph{A {B}urgess-like subconvex bound for
  twisted {$L$}-functions}, Forum Math. \textbf{19} (2007), no.~1, 61--106.

\bibitem{MR777279}
Joe~P. Buhler, Benedict~H. Gross, and Don~B. Zagier, \emph{On the conjecture of
  {B}irch and {S}winnerton-{D}yer for an elliptic curve of rank {$3$}}, Math.
  Comp. \textbf{44} (1985), no.~170, 473--481. \MR{MR777279 (86g:11037)}

\bibitem{CoMi}
J.~Cogdell and P.~Michel, \emph{On the complex moments of symmetric power
  {$L$}-functions at {$s=1$}}, Int. Math. Res. Not. (2004), no.~31, 1561--1617.
  \MR{MR2035301 (2005f:11094)}

\bibitem{DeIw}
J.-M. Deshouillers and H.~Iwaniec, \emph{Kloosterman sums and {F}ourier
  coefficients of cusp forms}, Invent. Math. \textbf{70} (1982/83), no.~2,
  219--288. \MR{MR684172 (84m:10015)}

\bibitem{Es}
T.~Estermann, \emph{On {K}loosterman's sum}, Mathematika \textbf{8} (1961),
  83--86. \MR{MR0126420 (23 \#A3716)}

\bibitem{55.0703.02}
Th. Estermann, \emph{{Vereinfachter Beweis eines Satzes von Kloosterman.}},
  Abhandlungen Hamburg \textbf{7} (1929), 82--98 (German).

\bibitem{GeJa}
Stephen Gelbart and Herv{\'e} Jacquet, \emph{A relation between automorphic
  representations of {${\rm GL}(2)$} and {${\rm GL}(3)$}}, Ann. Sci. \'Ecole
  Norm. Sup. (4) \textbf{11} (1978), no.~4, 471--542. \MR{MR533066 (81e:10025)}

\bibitem{Gu}
Ahmet~Muhtar G{\"u}lo{\u{g}}lu, \emph{Low-lying zeroes of symmetric power
  {$L$}-functions}, Int. Math. Res. Not. (2005), no.~9, 517--550.
  \MR{MR2131448}

\bibitem{MR1513069}
E.~Hecke, \emph{\"{U}ber die {B}estimmung {D}irichletscher {R}eihen durch ihre
  {F}unktionalgleichung}, Math. Ann. \textbf{112} (1936), no.~1, 664--699.
  \MR{MR1513069}

\bibitem{MR1513122}
\bysame, \emph{\"{U}ber {M}odulfunktionen und die {D}irichletschen {R}eihen mit
  {E}ulerscher {P}roduktentwicklung. {I}}, Math. Ann. \textbf{114} (1937),
  no.~1, 1--28. \MR{MR1513122}

\bibitem{MR1513142}
\bysame, \emph{\"{U}ber {M}odulfunktionen und die {D}irichletschen {R}eihen mit
  {E}ulerscher {P}roduktentwicklung. {II}}, Math. Ann. \textbf{114} (1937),
  no.~1, 316--351. \MR{MR1513142}

\bibitem{MR2166468}
C.~P. Hughes, \emph{Mock-{G}aussian behaviour}, Recent perspectives in random
  matrix theory and number theory, London Math. Soc. Lecture Note Ser., vol.
  322, Cambridge Univ. Press, Cambridge, 2005, pp.~337--355. \MR{MR2166468
  (2006g:11191)}

\bibitem{HuRu}
C.~P. Hughes and Z.~Rudnick, \emph{Linear statistics of low-lying zeros of
  {$L$}-functions}, Q. J. Math. \textbf{54} (2003), no.~3, 309--333.
  \MR{MR2013141 (2005a:11131)}

\bibitem{HuRu2}
\bysame, \emph{Mock-{G}aussian behaviour for linear statistics of classical
  compact groups}, J. Phys. A \textbf{36} (2003), no.~12, 2919--2932, Random
  matrix theory. \MR{MR1986399 (2004e:60012)}

\bibitem{Iw}
Henryk Iwaniec, \emph{Topics in classical automorphic forms}, Graduate Studies
  in Mathematics, vol.~17, American Mathematical Society, Providence, RI, 1997.
  \MR{MR1474964 (98e:11051)}

\bibitem{IwKo}
Henryk Iwaniec and Emmanuel Kowalski, \emph{Analytic number theory}, American
  Mathematical Society Colloquium Publications, vol.~53, American Mathematical
  Society, Providence, RI, 2004. \MR{MR2061214 (2005h:11005)}

\bibitem{IwLuSa}
Henryk Iwaniec, Wenzhi Luo, and Peter Sarnak, \emph{Low lying zeros of families
  of {$L$}-functions}, Inst. Hautes \'Etudes Sci. Publ. Math. (2000), no.~91,
  55--131 (2001). \MR{MR1828743 (2002h:11081)}

\bibitem{KaSa}
Nicholas~M. Katz and Peter Sarnak, \emph{Random matrices, {F}robenius
  eigenvalues, and monodromy}, American Mathematical Society Colloquium
  Publications, vol.~45, American Mathematical Society, Providence, RI, 1999.
  \MR{MR1659828 (2000b:11070)}

\bibitem{Ki}
Henry~H. Kim, \emph{Functoriality for the exterior square of {${\rm GL}\sb 4$}
  and the symmetric fourth of {${\rm GL}\sb 2$}}, J. Amer. Math. Soc.
  \textbf{16} (2003), no.~1, 139--183 (electronic), With appendix 1 by Dinakar
  Ramakrishnan and appendix 2 by Kim and Peter Sarnak. \MR{MR1937203
  (2003k:11083)}

\bibitem{KiSh2}
Henry~H. Kim and Freydoon Shahidi, \emph{Cuspidality of symmetric powers with
  applications}, Duke Math. J. \textbf{112} (2002), no.~1, 177--197.
  \MR{MR1890650 (2003a:11057)}

\bibitem{KiSh1}
\bysame, \emph{Functorial products for {${\rm GL}\sb 2\times{\rm GL}\sb 3$} and
  the symmetric cube for {${\rm GL}\sb 2$}}, Ann. of Math. (2) \textbf{155}
  (2002), no.~3, 837--893, With an appendix by Colin J. Bushnell and Guy
  Henniart. \MR{MR1923967 (2003m:11075)}

\bibitem{MR1915038}
E.~Kowalski, P.~Michel, and J.~VanderKam, \emph{Rankin-{S}elberg
  {$L$}-functions in the level aspect}, Duke Math. J. \textbf{114} (2002),
  no.~1, 123--191. \MR{MR1915038 (2004c:11070)}

\bibitem{MR848380}
Jean-Fran{\c{c}}ois Mestre, \emph{Courbes de {W}eil et courbes
  supersinguli\`eres}, Seminar on number theory, 1984--1985 (Talence,
  1984/1985), Univ. Bordeaux I, Talence, 1985, pp.~Exp.\ No.\ 23, 6.
  \MR{MR848380 (88b:11034)}

\bibitem{mil02a}
Stephen~D. Miller, \emph{The highest lowest zero and other applications of
  positivity}, Duke Math. J. \textbf{112} (2002), no.~1, 83--116. \MR{1 890
  648}

\bibitem{Mil}
Steven~J. Miller, \emph{One- and two-level densities for rational families of
  elliptic curves: evidence for the underlying group symmetries}, Compos. Math.
  \textbf{140} (2004), no.~4, 952--992. \MR{MR2059225 (2005c:11085)}

\bibitem{oma00}
Sami Omar, \emph{Majoration du premier z\'ero de la fonction z\^eta de
  {D}edekind}, Acta Arith. \textbf{95} (2000), no.~1, 61--65. \MR{MR1787205
  (2001h:11143)}

\bibitem{Sos00}
Alexander Soshnikov, \emph{The central limit theorem for local linear
  statistics in classical compact groups and related combinatorial identities},
  Ann. Probab. \textbf{28} (2000), no.~3, 1353--1370. \MR{MR1797877
  (2002f:15035)}

\bibitem{Sta97}
Richard~P. Stanley, \emph{Enumerative combinatorics. {V}ol. 1}, Cambridge
  Studies in Advanced Mathematics, vol.~49, Cambridge University Press,
  Cambridge, 1997, With a foreword by Gian-Carlo Rota, Corrected reprint of the
  1986 original. \MR{MR1442260 (98a:05001)}

\bibitem{Stan00}
\bysame, \emph{Enumerative combinatorics. {V}ol. 2}, Cambridge Studies in
  Advanced Mathematics, vol.~62, Cambridge University Press, Cambridge, 1999,
  With a foreword by Gian-Carlo Rota and appendix 1 by Sergey Fomin.
  \MR{MR1676282 (2000k:05026)}

\end{thebibliography}
%.......................................
\end{document}